\documentclass[a4paper,11pt,reqno]{amsart}

\newtheorem{theorem}[equation]{Theorem}
\newtheorem{lemma}[equation]{Lemma}
\newtheorem{corollary}[equation]{Corollary}
\newtheorem{proposition}[equation]{Proposition}

\numberwithin{equation}{section}

\theoremstyle{definition}

\newtheorem*{example*}{Example}
\newtheorem{example}[equation]{Example}

\newtheorem*{remark*}{Remark}

\newcommand{\bZ}{{\mathbb Z}}

\newcommand{\frg}{{\mathfrak g}}

\newcommand{\frt}{{\mathfrak t}}
\newcommand{\frs}{{\mathfrak s}}

\newcommand{\frv}{{\mathfrak v}}
\newcommand{\frl}{{\mathfrak l}}

\newcommand{\frn}{{\mathfrak n}}

\newcommand{\calT}{{\mathcal T}}

\newcommand{\calK}{{\mathcal K}}

\newcommand{\subo}{_{\bar 0}}
\newcommand{\subuno}{_{\bar 1}}

\providecommand{\espan}[1]{\text{span}\left\{ #1\right\}}

 \DeclareMathOperator{\tri}{\mathfrak{tri}}
  \DeclareMathOperator{\stri}{\mathfrak{stri}}
   \DeclareMathOperator{\lrt}{\mathfrak{lrt}}
    \DeclareMathOperator{\inlrt}{\mathfrak{inlrt}}
 \DeclareMathOperator{\fro}{\mathfrak{o}}
 \DeclareMathOperator{\frsl}{{\mathfrak{sl}}}

 \DeclareMathOperator{\frgl}{{\mathfrak{gl}}}

 \DeclareMathOperator{\frsu}{{\mathfrak{su}}}
 \DeclareMathOperator{\frstu}{{\mathfrak{stu}}}
 
 \DeclareMathOperator{\frd}{{\mathfrak{d}}}

 \DeclareMathOperator{\ad}{ad}
 \DeclareMathOperator{\der}{\mathfrak{der}}
 \DeclareMathOperator{\inder}{\mathfrak{inder}}

 \DeclareMathOperator{\Aut}{Aut}

\newenvironment{romanenumerate}
 {\begin{enumerate}
 \renewcommand{\itemsep}{4pt}
 }{\end{enumerate}}

\newenvironment{alphaenumerate}
 {\begin{enumerate}
 \renewcommand{\itemsep}{4pt}
 }{\end{enumerate}}

\begin{document}

\title[Lie algebras with $S_4$-action]{Lie algebras with $S_4$-action
and structurable algebras}

\author[Alberto Elduque]{Alberto Elduque$^{\star}$}
 \thanks{$^{\star}$ Supported by the Spanish Ministerio de
 Educaci\'{o}n y Ciencia
 and FEDER (MTM 2004-081159-C04-02) and by the
Diputaci\'on General de Arag\'on (Grupo de Investigaci\'on de
\'Algebra)}
 \address{Departamento de Matem\'aticas, Universidad de
Zaragoza, 50009 Zaragoza, Spain}
 \email{elduque@unizar.es}

\author[Susumu Okubo]{Susumu Okubo$^{\ast}$}
 \thanks{$^{\ast}$ Supported in part by U.S.~Department of Energy Grant No.
 DE-FG02-91ER40685.}
 \address{Department of Physics and Astronomy, University of
 Rochester, Rochester, NY 14627, USA}
 \email{okubo@pas.rochester.edu}

\date{\today}



\begin{abstract}
The normal symmetric triality algebras (STA's) and the normal Lie
related triple algebras (LRTA's) have been recently introduced by
the second author, in connection with the principle of triality. It
turns out that the unital normal LRTA's are precisely the
structurable algebras extensively studied by Allison.

It will be shown that the normal STA's (respectively LRTA's) are the
algebras that coordinatize those Lie algebras whose automorphism
group contains a copy of the alternating (resp. symmetric) group of
degree 4.
\end{abstract}

\maketitle


\section*{Introduction}

Over the years, many different constructions have been given of the
exceptional simple Lie algebras in Killing-Cartan's classification.
In 1966, Tits gave a unified construction of these algebras using a
couple of ingredients: a unital composition algebra and a simple
Jordan algebra of degree $3$ \cite{Tit66}. Even though the
construction is not symmetric, the outcome (the Magic Square)
presents a surprising symmetry.

A symmetric construction was obtained by Vinberg \cite{Vin66} (see
also \cite{OV94}) in terms of two unital composition algebras and
their Lie algebras of derivations. The two composition algebras play
the same role, and hence the symmetry of the construction. Vinberg's
construction was extended by Allison \cite{All91} in terms of
structurable algebras \cite{All78}, so that Vinberg's construction
becomes the particular case in which the structurable algebra is the
tensor product of two composition algebras. Moreover, Allison and
Faulkner \cite{AF93} gave a new version of Allison construction
which is based on three copies of a structurable algebra and a Lie
algebra of \emph{Lie related triples}. Quite recently, Barton and
Sudbery \cite{BS03} (see also the work by Landsberg and Manivel
\cite{LM02}, \cite{LM04}, and the survey by Baez \cite{Bae02}) gave
a simple recipe to obtain the Magic Square in terms of two unital
composition algebras and their triality Lie algebras, which in
perspective is subsumed in Allison-Faulkner's construction. Also,
since simpler formulas for triality are obtained by using the so
called \emph{symmetric composition algebras} instead of the
classical unital composition algebras, a version of Barton-Sudbery's
construction in terms of these latter algebras was given in
\cite{Eld04}, \cite{Eld05}. All these symmetric constructions of the
Magic Square provide $\bZ_2\times\bZ_2$-graded models of the Lie
algebras involved.

In \cite{Oku05}, the second author has studied two classes of
algebras: \emph{normal symmetric triality algebras} and \emph{normal
Lie related triple algebras}, whose defining conditions reflect
respectively the properties of the tensor products of two symmetric
composition algebras, needed in the construction of the Magic Square
in \cite{Eld04}, and the properties of structurable algebras needed
in Allison-Faulkner's construction \cite{AF93}. Algebras in both
classes are the building blocks of some $\bZ_2\times\bZ_2$-graded
Lie algebras.

The aim of this paper is to show that the algebras in these classes
are precisely those algebras that \emph{coordinatize} the Lie
algebras with an action of either the alternating or the symmetric
group of degree $4$. This possibility of relating these classes of
algebras with an action of these groups, as automorphisms of Lie
algebras, was suggested by the $S_4$-action on the exceptional
simple real Lie algebras considered by Loke \cite{Lok04}.

The paper is organized as follows. The first section deals with Lie
algebras with a subgroup of automorphisms isomorphic to the
alternating group of degree $4$. These Lie algebras are shown to be
coordinatized by the normal symmetric triality algebras. Then the
second section is devoted to Lie algebras with an action of the
larger symmetric group of degree $4$, which turn out to be
coordinatized by the normal Lie related triple algebras. The unital
such algebras are precisely the structurable algebras. Section 3
presents the main examples of algebras in these classes, which
include Jordan and structurable algebras, but also Lie algebras, Lie
triple systems and tensor products of symmetric composition
algebras. Finally, Section 4 shows how, under some restriction on
the ground field, Kantor's $5$-graded Lie algebras constructed from
structurable algebras are endowed with an action of $S_4$. This
gives a folding of the $5$-graded Lie algebras into very symmetric
$\bZ_2\times\bZ_2$-graded algebras, which were considered previously
by Allison and Faulkner \cite{AF93}.

\smallskip

\emph{Throughout the paper, all the algebras will be considered over
a ground field $F$ of characteristic $\ne 2$.}

\bigskip

\section{$A_4$-action and normal symmetric triality algebras}

Let $A_4$ denote the alternating group of degree $4$, and let $\frg$
be a Lie algebra endowed with a group homomorphism
\[
A_4\longrightarrow \Aut(\frg).
\]
The group $A_4$ is the semidirect product of the normal Klein's
$4$-group $V=\langle \tau_1,\tau_2\rangle$, where $\tau_1=(12)(34)$
and $\tau_2=(23)(14)$, and the cyclic group of order $3$,
$C_3=\langle \varphi\rangle$, where $\varphi=(123)$ ($1\mapsto
2\mapsto 3\mapsto 1$). Note that $\tau_1\tau_2=\tau_2\tau_1$,
$\varphi\tau_1=\tau_2\varphi$, and
$\varphi\tau_2=\tau_1\tau_2\varphi$. The same notation will be used
for the images of these elements in $\Aut(\frg)$.

The action of Klein's $4$-group gives a $\bZ_2\times\bZ_2$-grading
on $\frg$:
\begin{equation}\label{eq:z22grading}
\frg=\frt\oplus\frg_0\oplus\frg_1\oplus\frg_2,
\end{equation}
where
\begin{equation}\label{eq:tg0g1g2}
\begin{split}
\frt&=\{x\in\frg : \tau_1(x)=x,\, \tau_2(x)=x\}\
       (=\frg_{(\bar 0,\bar 0)}),\\
\frg_0&=\{x\in\frg : \tau_1(x)=x,\, \tau_2(x)=-x\}\
        (=\frg_{(\bar 1,\bar 0)}),\\
\frg_1&=\{x\in\frg : \tau_1(x)=-x,\, \tau_2(x)=x\}\
        (=\frg_{(\bar 0,\bar 1)}),\\
\frg_2&=\{x\in\frg : \tau_1(x)=-x,\, \tau_2(x)=-x\}\
        (=\frg_{(\bar 1,\bar 1)}).
\end{split}
\end{equation}
(Here, the subindices $0,1,2$ must be considered as the elements in
$\bZ_3$.)

Since $V$ is a normal subgroup of $A_4$, $\frt$ is invariant under
$\varphi$. Also, for any $x\in\frg_0$,
\[
\varphi(x)=\begin{cases}
\varphi\tau_1(x)=\tau_2\varphi(x),\\
-\varphi\tau_2(x)=-\tau_1\tau_2\varphi(x)=-\tau_1\varphi\tau_1(x)
   =-\tau_1\varphi(x),
\end{cases}
\]
so that $\varphi(\frg_0)\subseteq \frg_1$ and, in the same vein, one
gets $\varphi(\frg_i)\subseteq \frg_{i+1}$ for any $i$ (indices
modulo $3$).

Let $A=\frg_0$, and for any $x\in A$ consider the elements
\[
\iota_0(x)=x\in \frg_0,\quad \iota_1(x)=\varphi\iota_0(x)\in
\frg_1,\quad \iota_2(x)=\varphi^2\iota_0(x)\in\frg_2.
\]
Thus
\[
\frg=\frt\oplus \bigl(\oplus_{i=0}^2\iota_i(A)\bigr).
\]

Since \eqref{eq:z22grading} is a grading over $\bZ_2\times\bZ_2$,
 a bilinear multiplication $*$ can be defined on $A$ by means of
\begin{equation}\label{eq:iota1iota2}
[\iota_1(x),\iota_2(y)]=\iota_0(x*y),
\end{equation}
for any $x,y\in A$. Also, the fact that $\varphi$ is an automorphism
shows that \eqref{eq:iota1iota2} is equivalent to
\begin{equation}\label{eq:iotaiiplus1}
[\iota_i(x),\iota_{i+1}(y)]=\iota_{i+2}(x*y),
\end{equation}
for any $x,y\in A$ and $i\in\bZ_3$.

\medskip

Given any algebra $(A,*)$, consider its \emph{symmetric triality Lie
algebra} (see \cite{Oku05}):
\begin{multline*}
\stri(A,*)=\{ (d_0,d_1,d_2)\in\frgl(A)^3:
d_i(x*y)=d_{i+1}(x)*y+x*d_{i+2}(y)\\
\text{for any }i=0,1,2\text{ and }x,y\in A\}.
\end{multline*}
This is a Lie subalgebra of $\frgl(A)^3$ (with componentwise Lie
bracket), where $\frgl(A)$ is the Lie algebra of endomorphisms of
the vector space $A$.

\medskip

Let $A=\frg_0$ as above. Three representations of $\frt$ on $A$:
$\rho_i:\frt\rightarrow \frgl(A)$ ($i=0,1,2$), are obtained by means
of:
\[
\iota_i\bigl(\rho_i(d)(x)\bigr)=[d,\iota_i(x)],
\]
for any $d\in\frt$ and $x\in A$. Note that
\[
\begin{split}
\iota_{i+1}\Bigl(\rho_{i+1}\bigl(\varphi(d)\bigr)(x)\Bigr)&=
 \bigl[\varphi(d),\iota_{i+1}(x)\bigr]=\bigl[\varphi(d),\varphi(\iota_i(x))\bigr]\\
 &=\varphi\Bigl(\bigl[d,\iota_i(x)\bigr]\Bigr)=
  \varphi\Bigl(\iota_i\bigl(\rho_i(d)(x)\bigr)\Bigr)\\
 &=\iota_{i+1}\bigl(\rho_i(d)(x)\bigr),
\end{split}
\]
so $\rho_{i+1}\varphi=\rho_i$, or
\begin{equation}\label{eq:rhovarphi}
\rho_i\varphi^j=\rho_{i-j},
\end{equation}
for any $i,j\in \bZ_3$.

\smallskip

Putting together $\rho_0$, $\rho_1$ and $\rho_2$, there appears a
Lie algebra homomorphism:
\[
\begin{split}
\rho: \frt&\longrightarrow \frgl(A)^3\\
d&\mapsto \bigl(\rho_0(d),\rho_1(d),\rho_2(d)\bigr).
\end{split}
\]
Then:

\begin{proposition}\label{pr:rhofrt}
Under the hypotheses above, $\rho(\frt)$ is contained in
$\stri(A,*)$.
\end{proposition}
\begin{proof}
For any $d\in\frt$, $x,y\in A$ and $i\in\bZ_3$,
\[
\begin{split}
\iota_i\bigl(\rho_i(d)(x*y)\bigr)&=
 \bigl[d,\iota_i(x*y)\bigr]\\
 &=\bigl[d,\bigl[\iota_{i+1}(x),\iota_{i+2}(x)\bigr]\bigr]\\
 &=\bigl[\bigl[d,\iota_{i+1}(x)\bigr],\iota_{i+2}(y)\bigr]
   +\bigl[\iota_{i+1}(x),\bigr[d,\iota_{i+2}(y)\bigr]\bigr]\\
 &=\bigl[\iota_{i+1}\bigl(\rho_{i+1}(d)(x)\bigr),\iota_{i+2}(y)\bigr]
   +\bigl[\iota_{i+1}(x),\iota_{i+2}\bigl(\rho_{i+2}(d)(y)\bigr)\bigr]\\
 &=\iota_i\Bigl(\rho_{i+1}(d)(x)*y+x*\rho_{i+2}(d)(y)\Bigr),
\end{split}
\]
which proves the result.
\end{proof}

Therefore, $\rho$ becomes a homomorphism of Lie algebras:
\[
\rho: \frt\longrightarrow \stri(A,*).
\]
Note that $\stri(A,*)$ has a natural order $3$ automorphism
\begin{equation}\label{eq:theta}
\theta: (d_0,d_1,d_2)\mapsto (d_2,d_0,d_1)
\end{equation}
and
\eqref{eq:rhovarphi} is equivalent to:
\begin{equation}\label{eq:rhovarphitheta}
\rho\varphi=\theta\rho.
\end{equation}

Consider now the skew-symmetric linear map
\[
\begin{split}
\delta: A\times A&\longrightarrow \stri(A,*)\\
(x,y)\ &\mapsto \delta(x,y)=\rho\bigl([\iota_0(x),\iota_0(y)]\bigr),
\end{split}
\]
and denote by $\delta_i(x,y)$ the $i^{\text{th}}$ component of
$\delta(x,y)$
($\delta_i(x,y)=\rho_i\bigl([\iota_0(x),\iota_0(y)\bigr]$), so that
$\delta(x,y)=\bigl(\delta_0(x,y),\delta_1(x,y),\delta_2(x,y)\bigr)$.
Note that $[\frg_i,\frg_i]\subseteq \frt$ ($i=0,1,2$).

\begin{theorem}\label{th:deltas}
Under the conditions above, for any $a,b,x,y,z\in A$ and
$i,j\in\bZ_3$:
\begin{romanenumerate}
\item
$\bigl[\delta_i(a,b),\delta_j(x,y)\bigr]
 =\delta_j\bigl(\delta_{i-j}(a,b)(x),y\bigr)+
 \delta_j(x,\delta_{i-j}(a,b)(y)\bigr)$,
\item $\delta_0(x,y*z)+\delta_{-1}(y,z*x)+\delta_{-2}(z,x*y)=0$,
\item $\delta_0(x,y)(z)+\delta_0(y,z)(x)+\delta_0(z,x)(y)=0$,
\item $\delta_1(x,y)=r_yl_x-r_xl_y$ ($l_x: z\mapsto x*z$,
 $r_x:z\mapsto z*x$),
\item $\delta_2(x,y)=l_yr_x-l_xr_y$.
\end{romanenumerate}
\end{theorem}

\begin{proof}
For any $d\in \frt$ and $x,y\in A$:
\begin{equation}\label{eq:rhodxy}
\begin{split}
\bigl[\rho(d),\delta(x,y)\bigr]&=
 \rho\Bigl(\bigl[d,\bigl[\iota_0(x),\iota_0(y)\bigr]\bigr]\Bigr)\\
 &=\rho\Bigl(\bigl[\bigl[d,\iota_0(x)\bigr],\iota_0(y)\bigr]
   +\bigl[\iota_0(x),\bigl[d,\iota_0(y)\bigr]\bigr]\Bigr)\\
 &=\rho\Bigl(\bigl[\iota_0\bigl(\rho_0(d)(x)\bigr),\iota_0(y)\bigr]
   +\bigl[\iota_0(x),\iota_0\bigl(\rho_0(d)(y)\bigr)\bigr]\Bigr)\\
 &=\delta\bigl(\rho_0(d)(x),y\bigr)+\delta\bigl(x,\rho_0(d)(y)\bigr).
\end{split}
\end{equation}
Now, for $d=[\iota_i(a),\iota_i(b)]$, using
\eqref{eq:rhovarphitheta} one gets
\[
\begin{split}
\rho(d)&=\rho\varphi^i\bigl([\iota_0(a),\iota_0(b)]\bigr)\\
 &=\theta^i\rho\bigl([\iota_0(a),\iota_0(b)]\bigr)\\
 &=\theta^i\bigl(\delta(a,b)\bigr)\\
 &=\bigl(\delta_{-i}(a,b),\delta_{1-i}(a,b),\delta_{2-i}(a,b)\bigr),
\end{split}
\]
and the $j^{\text{th}}$ component of \eqref{eq:rhodxy} becomes
\[
\bigl[\delta_{j-i}(a,b),\delta_j(x,y)\bigr]=
 \delta_j\bigl(\delta_{-i}(a,b)(x),y\bigr)
 +\delta_j\bigl(x,\delta_{-i}(a,b)(y)\bigr),
\]
which is equivalent to the assertion in item (i).

Now, the Jacobi identity implies:
\[
\begin{split}
0&=\bigl[\iota_0(x),[\iota_1(y),\iota_2(z)]\bigr] +
  \bigl[\iota_1(y),[\iota_2(z),\iota_0(x)]\bigr] +
  \bigl[\iota_2(z),[\iota_0(x),\iota_1(y)]\bigr]\\
 &=[\iota_0(x),\iota_0(y*z)]+[\iota_1(y),\iota_1(z*x)]
   +[\iota_2(z),\iota_2(x*y)]\\
 &=[\iota_0(x),\iota_0(y*z)]+\varphi\bigl([\iota_0(y),\iota_0(z*x)]\bigr)
   +\varphi^2\bigl([\iota_0(z),\iota_0(x*y)]\bigr).
\end{split}
\]
Apply $\rho$ and use \eqref{eq:rhovarphitheta} to obtain:
\[
\delta(x,y*z)+\theta\bigl(\delta(y,z*x)\bigr)+
 \theta^2\bigl(\delta(z,x*y)\bigr)=0,
\]
whose first component gives (ii).

Also,
\[
\begin{split}
0&=\bigl[[\iota_0(x),\iota_0(y)],\iota_0(z)\bigr]+
  \bigl[[\iota_0(y),\iota_0(z)],\iota_0(x)\bigr]+
  \bigl[[\iota_0(z),\iota_0(x)],\iota_0(y)\bigr]\\
   &=\iota_0\bigl(\delta_0(x,y)(z)+\delta_0(y,z)(x)+\delta_0(z,x)(y)\bigr),
\end{split}
\]
whence (iii),
\[
\begin{split}
\iota_1\bigl(\delta_1(x,y)(z)\bigr)&=
  \bigl[[\iota_0(x),\iota_0(y)],\iota_1(z)\bigr]\\
  &=\bigl[[\iota_0(x),\iota_1(z)],\iota_0(y)\bigr]+
     \bigl[\iota_0(x),[\iota_0(y),\iota_1(z)]\bigr]\\
  &=\bigl[\iota_2(x*z),\iota_0(y)\bigr]
     +\bigl[\iota_0(x),\iota_2(y*z)\bigr]\\
  &=\iota_1\Bigl((x*z)*y-(y*z)*x\Bigr),
\end{split}
\]
hence (iv), and
\[
\begin{split}
\iota_2\bigl(\delta_2(x,y)(z)\bigr)&=
   \bigl[[\iota_0(x),\iota_0(y)],\iota_2(z)\bigr]\\
    &=\bigl[[\iota_0(x),\iota_2(z)],\iota_0(y)\bigr]
     +\bigl[\iota_0(x),[\iota_0(y),\iota_2(z)]\bigr]\\
    &=-\bigl[\iota_1(z*x),\iota_0(y)\bigr]
      -\bigl[\iota_0(x),\iota_1(z*y)\bigr]\\
    &=\iota_2\Bigl(y*(z*x)-x*(z*y)\Bigr),
\end{split}
\]
which proves (v).
\end{proof}

\begin{remark*} Note that if $A*A=A$, then $\delta_0$ is determined
by $\delta_1$ and $\delta_2$ because of item (ii) in Theorem
\ref{th:deltas} (or by ``triality'':
$\delta_0(x,y)(u*v)=\delta_1(x,y)(u)*v+u*\delta_2(x,y)(v)$). Also,
if $\{x\in A:A*x=0\}=0$, $\delta_0$ is determined by $\delta_1$ and
$\delta_2$ by triality, since
$\delta_1(x,y)(u*v)=\delta_2(x,y)(u)*v+u*\delta_0(x,y)(v)$, and the
same happens if $\{x\in A: x*A=0\}=0$.
\end{remark*}

\medskip

Conditions (i)--(v) in Theorem \ref{th:deltas} are precisely the
conditions (2.15) in \cite{Oku05} defining a \emph{normal symmetric
triality algebra}. Therefore:

\begin{corollary}
Under the conditions of Theorem \ref{th:deltas}, the algebra $(A,*)$
is a normal symmetric triality  algebra.
\end{corollary}

\begin{corollary}
Let $(A,*)$ be a normal symmetric triality algebra with respect to
the skew-symmetric bilinear map $\delta:A\times A\rightarrow
\stri(A,*)$. Then $(A,*)$ satisfies the degree $5$ identity:
\begin{equation}\label{eq:degree5}
\begin{split}
0&=\bigl((x*u)*(y*z)\bigr)*v-\bigl(((y*z)*u)*x\bigr)*v\\
 &\quad +u*\bigl((y*z)*(v*x)\bigr)-u*\bigl(x*(v*(y*z))\bigr)\\
 &\quad +(z*x)*\bigl((u*v)*y\bigr)-y*\bigl((u*v)*(z*x)\bigr)\\
 &\quad +\bigl(z*(u*v)\bigr)*(x*y)-\bigl((x*y)*(u*v)\bigr)*z,
\end{split}
\end{equation}
for any $u,v,x,y,z\in A$.
\end{corollary}
\begin{proof}
By `triality', for any $u,v,x,y,z\in A$:
\[
\delta_0(x,y*z)(u*v)=\delta_1(x,y*z)(u)*v+u*\delta_2(x,y*z)(v),
\]
while item (ii) in Theorem \ref{th:deltas} gives:
\[
\delta_0(x,y*z)(u*v)=-\delta_2(y,z*x)(u*v)-\delta_1(z,x*y)(u*v).
\]
Hence,
\begin{multline*}
0=\delta_1(x,y*z)(u)*v+u*\delta_2(x,y*z)(v)\\
 +\delta_2(y,z*x)(u*v)+\delta_1(z,x*y)(u*v).
\end{multline*}
Expanding this last equation, by means of items (iv) and (v) of
Theorem \ref{th:deltas}, gives \eqref{eq:degree5}.
\end{proof}

\medskip

The computations in the proof of Theorem \ref{th:deltas} can be
reversed to get a sort of converse. The straightforward proof is
omitted.

\begin{theorem}\label{th:reverseSTA}
Let $(A,*)$ be a nonzero normal symmetric triality algebra with
respect to the skew-symmetric bilinear map $\delta: A\times
A\rightarrow \stri(A,*)$. Then
$\frt=\sum_{i=0}^2\theta^i\bigl(\delta(A,A)\bigr)$ is a Lie
subalgebra of $\stri(A,*)$, and $\theta^i\bigl(\delta(A,A)\bigr)$ is
an ideal of $\frt$ for any $i=0,1,2$. Moreover, consider three
copies $\iota_i(A)$ of $A$ ($i=0,1,2$) and define an anticommutative
multiplication on
\begin{equation}\label{eq:gastar}
\frg(A,*)=\frt\oplus\bigl(\oplus_{i=0}^2 \iota_i(A)\bigr)
\end{equation}
by means of
\begin{itemize}
\renewcommand*{\itemsep}{4pt}
\item $\frt$ is a subalgebra of $\frg(A,*)$,
\item
$\bigl[(d_0,d_1,d_2),\iota_i(x)\bigr]=\iota_i\bigl(d_i(x)\bigr)$,
for any $(d_0,d_1,d_2)\in \frt$, $x\in A$ and $i\in\bZ_3$,
\item $\bigl[\iota_i(x),\iota_{i+1}(y)]=\iota_{i+2}(x*y)$ for any
$x,y\in A$ and $i\in\bZ_3$,
\item
$\bigl[\iota_i(x),\iota_i(y)\bigr]=\theta^i\bigl(\delta(x,y)\bigr)$,
for any $x,y\in A$ and $i\in\bZ_3$.
\end{itemize}
Then $\frg(A,*)$ is a Lie algebra and the alternating group $A_4$
embeds as a subgroup of $\Aut\frg(A,*)$ by means of
\begin{itemize}
\renewcommand*{\itemsep}{4pt}
\item The restrictions of the elements of the $4$-group $V$ to $\frt$ are trivial:
 $\tau_1\vert_\frt = id=\tau_2\vert_\frt$. Moreover,
 $\varphi\vert_\frt=\theta\vert_\frt$.
\item $\varphi\bigl(\iota_i(x)\bigr)=\iota_{i+1}(x)$ for any $x\in
A$ and $i\in\bZ_3$.
\item For any $x\in A$, $\tau_1\bigl(\iota_0(x)\bigr)=\iota_0(x)$
and $\tau_1\bigl(\iota_i(x)\bigr)=-\iota_i(x)$ for $i=1,2$, while
$\tau_2\bigl(\iota_1(x)\bigr)=\iota_1(x)$ and
$\tau_2\bigl(\iota_i(x)\bigr)=-\iota_i(x)$ for $i=0,2$.
\end{itemize}
\end{theorem}

\medskip

We finish the section with a result relating the simplicity of a
normal symmetric triality algebra and the associated Lie algebra:

\begin{theorem}\label{th:simpleSTA}
Let $(A,*)$ be a normal symmetric triality algebra with $A*A\ne 0$.
Then $(A,*)$ is simple if and only if the Lie algebra $\frg(A,*)$
defined in \eqref{eq:gastar} is simple as a Lie algebra with
$A_4$-action (that is, it does not contain any proper ideal
invariant under the action of $A_4$).
\end{theorem}
\begin{proof}
Let us show first that if $\frg(A,*)$ is simple as a Lie algebra
with $A_4$-action, then $A=A*A$. To see this, note that the subspace
\[
 \Bigl(\sum_{i=0}^2[\iota_i(A*A),\iota_i(A)]\Bigr) \oplus
 \Bigl(\oplus_{i=0}^2\iota_i(A*A)\Bigr)
\]
is invariant under the action of $A_4$, and it is closed under the
adjoint action of $\iota_i(A)$, $i=0,1,2$, which generate
$\frg(A,*)$. Actually,
\[
\bigl[[\iota_i(A*A),\iota_i(A)],\iota_{i\pm 1}(A)\bigr]\subseteq
  [\iota_{i\mp 1}(A*A),\iota_i(A)]\subseteq \iota_{i\pm 1}(A*A),
\]
because of the Jacobi identity, and also
\[
\bigl[[\iota_i(A*A),\iota_i(A)],\iota_i(A)\bigr]\subseteq
\iota_i\bigl(\delta_0(A*A,A)(A)\bigr)\subseteq \iota_i(A*A),
\]
because $\delta_0(A*A,A)\subseteq \delta_1(A*A,A)+\delta_2(A*A,A)$
by Theorem \ref{th:deltas}. Hence, the subspace above is an ideal of
$\frg(A,*)$ invariant under the action of $A_4$, so it is the whole
$\frg(A,*)$, and this shows that $A*A=A$.

Now, since $\delta_0(A,A)=\delta_0(A*A,A)\subseteq
\delta_1(A,A)+\delta_2(A,A)$, it is contained in the Lie
multiplication algebra of $A$ (the Lie subalgebra of $\frgl(A)$
generated by the left and right multiplications) by Theorem
\ref{th:deltas}.

Assume that $\frg(A,*)$ is simple as an algebra with $A_4$-action
and let $0\ne I$ be an ideal of $(A,*)$. By the above,
$\delta_0(A,I)(A)\subseteq I$. Then
\[
\Bigl(\sum_{i=0}^2[\iota_i(A),\iota_i(I)]\Bigr)\oplus
\iota_0(I)\oplus\iota_1(I)\oplus\iota_2(I)
\]
is closed under the adjoint action of $\iota_i(A)$ for $i=0,1,2$
since, for $i\ne j\ne k\ne i$,
\[
\begin{split}
\bigl[\iota_i(I),\iota_j(A)\bigr]&\subseteq
\iota_k(I*A+A*I)\subseteq \iota_k(I),\\
\bigl[[\iota_i(A),\iota_i(I)],\iota_j(A)\bigr]
 &\subseteq \bigl[\iota_k(A),\iota_i(I)\bigr]+
   \bigl[\iota_i(A),\iota_k(I)\bigr]\subseteq \iota_j(I),\\
 \bigl[[\iota_i(A),\iota_i(I)],\iota_i(A)\bigr]&\subseteq
 \iota_i\bigl(\delta_0(A,I)(A)\bigr)\subseteq \iota_i(I).
\end{split}
\]
Thus, since the $\iota_i(A)$'s generate $\frg(A,*)$,  this is a
nonzero ideal and hence equals the whole $\frg(A,*)$. Therefore,
$I=A$.

Finally, assume that $(A,*)$ is simple and let
$\frn=(\frn\cap\frt)\oplus
\Bigl(\oplus_{i=1}^2\bigl(\frn\cap\iota_i(A)\bigr)\Bigr)$ be an
ideal of $\frg(A,*)$ invariant under the action of $A_4$. For any
$i$, let $I_i=\{x\in A: \iota_i(x)\in \frn\}$. Then for any
$i\in\bZ_3$,
\[
\begin{split}
\bigl[\frn\cap\iota_i(A),\iota_{i+1}(A)\bigr]&\subseteq
  \frn\cap\iota_{i+2}(A),\\
\bigl[\iota_i(A),\frn\cap\iota_{i+1}(A)\bigr]
  &\subseteq \frn\cap\iota_{i+2}(A),
\end{split}
\]
which implies that $I_i*A\subseteq I_{i+2}$ and $A*I_{i+1}\subseteq
I_{i+2}$. But $\frn$ is invariant under the automorphism $\varphi$,
and this shows that $I_1=I_2=I_3$ is an ideal of $A$. If this is the
whole $A$, $\frn=\frg(A,*)$, while if this ideal is $0$,
$\frn\subseteq \frt$, which acts faithfully on
$\oplus_{i=0}^2\iota_i(A)$. However,
$\bigl[\frn,\iota_i(A)\bigr]\subseteq
\frn\cap\iota_i(A)=\iota_i(I_i)=0$. Hence $\frn=0$.
\end{proof}

The restriction $A*A\ne 0$ in Theorem \ref{th:simpleSTA} is
necessary, as shown by Example \ref{ex:LTS}.

\begin{remark*} Let $G$ be a group of automorphisms of a finite
dimensional Lie algebra $L$, and assume that $L$ is simple as a Lie
algebra with $G$-action. Then, if $I$ is a minimal ideal of $L$,
then $\sum_{\sigma\in G}\sigma(I)$ is an ideal, invariant under the
action of $G$. Hence $L=\sum_{\sigma\in G}\sigma(I)$, so $L$ is
completely reducible as a module over itself (the adjoint module).
It follows that $L=\sigma_1(I)\oplus\cdots\oplus\sigma_r(I)$ for
some $\sigma_1=1,\sigma_2,\ldots,\sigma_r\in G$, and each
$\sigma_j(I)$ is simple (by minimality of $I$). Therefore $L$ is a
direct sum of simple ideals.

In particular, this applies to the Lie algebras $\frg(A,*)$ for
simple $(A,*)$ in Theorem \ref{th:simpleSTA}.
\end{remark*}

\bigskip

\section{$S_4$-action and normal Lie related triple algebras}

In this section $\frg$ will be a Lie algebra endowed with a group
homomorphism
\[
S_4\longrightarrow \Aut(\frg),
\]
where $S_4$ is the symmetric group of degree $4$, which is the
semidirect product of $A_4$ and the cyclic subgroup of order $2$
generated by the transposition $\tau=(12)$. Besides,
$\tau_1\tau=\tau\tau_1$, $\tau_2\tau=\tau\tau_2\tau_1$, and
$\tau\varphi=\varphi^2\tau$. Because of the results in the previous
Section,
\[
\frg=\frt\oplus\bigl(\oplus_{i=0}^2 \iota_i(A)\bigr),
\]
where $(A,*)$ is a normal symmetric triality algebra (normal STA for
short).

Since Klein's $4$-group $V$ is a normal subgroup of $S_4$, $\frt$ is
invariant under $\tau$, and hence under $S_4$. For any $x\in
\frg_0=\iota_0(A)=\{x\in\frg: \tau_1(x)=x,\, \tau_2(x)=-x\}$,
\[
\tau_1\tau(x)=\tau\tau_1(x)=\tau(x),
\]
while
\[
\tau_2\tau(x)=\tau\tau_2\tau_1(x)=\tau\tau_2(x)=-\tau(x).
\]
Hence $\tau(\frg_0)\subseteq\frg_0$. Also, for any
$x\in\frg_1=\iota_1(A)=\{ x\in\frg: \tau_1(x)=-x,\,\tau_2(x)=x\}$,
\[
\begin{split}
\tau_1\tau(x)&=\tau\tau_1(x)=-\tau(x),\\
\tau_2\tau(x)&=\tau\tau_2\tau_1(x)=-\tau\tau_2(x)=-\tau(x).
\end{split}
\]
Therefore, $\tau(\frg_1)\subseteq \frg_2$, and also
$\tau(\frg_2)\subseteq \frg_1$.

An involutive linear map $A\rightarrow A$, $x\mapsto \bar x$, can be
defined  by means of
\[
\tau\bigl(\iota_0(x)\bigr)=-\iota_0(\bar x).
\]

\begin{proposition}\label{pr:bar}
The map $x\mapsto \bar x$ is an involution of $(A,*)$.
\end{proposition}
\begin{proof} First, we have $\bar{\bar x}=x$, since we calculate
\[
\iota_0(x)=\tau^2\bigl(\iota_0(x)\bigr)
 =-\tau\bigl(\iota_0(\bar x)\bigr)=\iota_0(\bar{\bar x}).
\]
Next, since $\tau\varphi=\varphi^2\tau$, for any $x\in A$,
\[
\begin{split}
\tau\iota_1(x)&=\tau\varphi\bigl(\iota_0(x)\bigr)
 =\varphi^2\tau\bigl(\iota_0(x)\bigr)=
 -\varphi^2\bigl(\iota_0(\bar x)\bigr)=-\iota_2(\bar x),\\
\tau\iota_2(x)&=\tau\varphi^2\bigl(\iota_0(x)\bigr)
 =\varphi\tau\bigl(\iota_0(x)\bigr)=
 -\varphi\bigl(\iota_0(\bar x)\bigr)=-\iota_1(\bar x).
\end{split}
\]
Thus, for any $x,y\in A$, apply the automorphism $\tau$ to
$\bigl[\iota_0(x),\iota_1(y)\bigr]=\iota_2(x*y)$ to get
$\bigl[\iota_0(\bar x),\iota_2(\bar
y)\bigr]=-\iota_1(\overline{x*y})$, or $-\iota_1(\bar y*\bar
x)=-\iota_1(\overline{x*y})$. Hence $\overline{x*y}=\bar y*\bar x$,
as required.
\end{proof}

Define a new multiplication on $A$ by means of:
\[
x\cdot y=\overline{x*y}=\bar y*\bar x,
\]
for any $x,y\in A$. Then $x\mapsto \bar x$ is an involution too of
$(A,\cdot)$. The second author has shown \cite[(1.16)]{Oku05} that
\[
\stri(A,*)=\lrt(A,\cdot,\bar\ ),
\]
where
\begin{multline*}
\lrt(A,\cdot,\bar\ )=\{(d_0,d_1,d_2)\in \frgl(A)^3 : \bar d_i(x\cdot
y)
  =d_{i+1}(x)\cdot y+x\cdot d_{i+2}(y)\\
  \text{for any }x,y\in A\text{ and }i\in\bZ_3\},
\end{multline*}
with $\bar d(x)=\overline{d(\bar x)}$ for any $d\in \frgl(A)$ and
$x\in A$.

The skew-symmetric bilinear map $\delta$ can be considered now as a
map
\[
\delta: A\times A\longrightarrow \lrt(A,\cdot,\bar\ ).
\]
As before, $\lrt(A,\cdot,\bar\ )=\stri(A,*)$ has the natural order
$3$ automorphism $\theta$  (see \eqref{eq:theta}) given by
\[
\theta(d_0,d_1,d_2)=(d_2,d_0,d_1),
\]
and also the order $2$ automorphism $\xi$  given by
\begin{equation}\label{eq:xi}
\xi(d_0,d_1,d_2)=(\bar d_0,\bar d_2,\bar d_1),
\end{equation}
which satisfies,
\begin{equation}\label{eq:rhotauxi}
\rho\tau=\xi\rho.
\end{equation}
This is shown with the same sort of arguments leading to
\eqref{eq:rhovarphi}.

\medskip

Denote by $L_x$ and $R_x$ the left and right multiplications by an
element $x$ in $(A,\cdot)$. Hence $L_x(y)=x\cdot
y=\overline{x*y}=\bar y*\bar x$, so $L_x=\nu l_x=r_{\bar x}\nu$,
where $\nu(x)=\bar x$. Also $R_x=\nu r_x=l_{\bar x}\nu$.

The maps $\delta_1(x,y)$, $\delta_2(x,y)$ in Theorem \ref{th:deltas}
become now:
\[
\left\{\begin{aligned}
 \delta_1(x,y)&=r_yl_x-r_xl_y=L_{\bar y}L_x-L_{\bar x}L_y,\\
 \delta_2(x,y)&=l_yr_x-l_xr_y=R_{\bar y}R_x-R_{\bar x}R_y.
 \end{aligned}\right.
\]
Since $\overline{R_x}=L_{\bar x}$ and $\overline{L_x}=R_{\bar x}$ in
any algebra with involution, it is clear that
$\overline{\delta_1(x,y)}=\delta_2(\bar x,\bar y)$ for any $x,y$.

Also, apply the automorphism $\tau$ to
$\bigl[[\iota_0(x),\iota_0(y)],\iota_0(z)\bigr]=
\iota_0\bigl(\delta_0(x,y)(z)\bigr)$ to obtain
\[
-\bigl[[\iota_0(\bar x),\iota_0(\bar y)],\iota_0(\bar z)\bigr]
 =-\iota_0\bigl(\overline{\delta_0(x,y)(z)}\bigr),
\]
so $\overline{\delta_0(x,y)}=\delta_0(\bar x,\bar y)$.

Hence, Theorem \ref{th:deltas} immediately implies the following:

\begin{theorem}\label{th:deltasNLRTA}
Under the conditions above, for any $a,b,x,y,z\in A$ and
$i,j\in\bZ_3$:
\begin{romanenumerate}
\item
$\bigl[\delta_i(a,b),\delta_j(x,y)\bigr]
 =\delta_j\bigl(\delta_{i-j}(a,b)(x),y\bigr)+
 \delta_j(x,\delta_{i-j}(a,b)(y)\bigr)$,
\item $\delta_0(\bar x,y\cdot z)+\delta_{1}(\bar y,z\cdot x)
   +\delta_{2}(\bar z,x\cdot y)=0$,
\item $\delta_0(x,y)(z)+\delta_0(y,z)(x)+\delta_0(z,x)(y)=0$,
\item $\delta_1(x,y)=L_{\bar y}L_x-L_{\bar x}L_y$,
\item $\delta_2(x,y)=R_{\bar y}R_x-R_{\bar x}R_y$,
\item $\overline{\delta_i(x,y)}=\delta_{-i}(\bar x,\bar y)$ (or
$\xi\bigl(\delta(x,y)\bigr)=\delta(\bar x,\bar y)$).
\end{romanenumerate}
\end{theorem}

Conditions (i)--(vi) above are precisely the conditions (2.34) in
\cite{Oku05} defining a \emph{normal Lie related triple algebra} (or
normal LRTA for short). Therefore:

\begin{corollary}
Under the hypothesis above, $(A,\cdot,\bar{\ })$ is a normal LRTA.
\end{corollary}

And, as in Section 2, everything can be reversed to get:

\begin{theorem}\label{th:reverseLRTA}
Let $(A,\cdot,\bar{\ })$ be a nonzero normal LRTA with respect to
the skew-symmetric bilinear map $\delta: A\times A\rightarrow
\lrt(A,\cdot,\bar{\ })$. Then
$\frt=\sum_{i=0}^2\theta^i\bigl(\delta(A,A)\bigr)$ is a Lie
subalgebra of $\lrt(A,\cdot,\bar\ )$, and
$\theta^i\bigl(\delta(A,A)\bigr)$ is an ideal of $\frt$ for any
$i=0,1,2$. Moreover, consider three copies $\iota_i(A)$ of $A$
($i=0,1,2$) and define an anticommutative multiplication on
\[
\frg(A,\cdot,\bar\ )=\frt\oplus\bigl(\oplus_{i=0}^2 \iota_i(A)\bigr)
\]
by means of
\begin{itemize}
\renewcommand*{\itemsep}{4pt}
\item $\frt$ is a subalgebra of $\frg(A,\cdot,\bar\ )$,
\item
$\bigl[(d_0,d_1,d_2),\iota_i(x)\bigr]=\iota_i\bigl(d_i(x)\bigr)$,
for any $(d_0,d_1,d_2)\in \frt$, $x\in A$ and $i\in\bZ_3$,
\item $\bigl[\iota_i(x),\iota_{i+1}(y)]=\iota_{i+2}(\overline{x\cdot y})$ for any
$x,y\in A$ and $i\in\bZ_3$,
\item
$\bigl[\iota_i(x),\iota_i(y)\bigr]=\theta^i\bigl(\delta(x,y)\bigr)$,
for any $x,y\in A$ and $i\in\bZ_3$.
\end{itemize}
Then $\frg(A,\cdot,\bar\ )$ is a Lie algebra and the symmetric group
$S_4$ embeds as a subgroup of $\Aut\frg(A,\cdot)$ by means of
\begin{itemize}
\renewcommand*{\itemsep}{4pt}
\item The restrictions of the elements of the $4$-group $V$ to $\frt$ are trivial:
 $\tau_1\vert_\frt = id=\tau_2\vert_\frt$. Moreover,
 $\varphi\vert_\frt=\theta\vert_\frt$ and
 $\tau\vert_\frt=\xi\vert_\frt$. ($\theta$ and $\xi$ as in
 \eqref{eq:theta} and \eqref{eq:xi}.)
\item $\varphi\bigl(\iota_i(x)\bigr)=\iota_{i+1}(x)$ for any $x\in
A$ and $i\in\bZ_3$.
\item For any $x\in A$, $\tau_1\bigl(\iota_0(x)\bigr)=\iota_0(x)$,
$\tau_1\bigl(\iota_i(x)\bigr)=-\iota_i(x)$ for $i=1,2$, while
$\tau_2\bigl(\iota_1(x)\bigr)=\iota_1(x)$,
$\tau_2\bigl(\iota_i(x)\bigr)=-\iota_i(x)$ for $i=0,2$.
\item For any $x\in A$, $\tau\bigl(\iota_0(x)\bigr)=-\iota_0(\bar
x)$, $\tau\bigl(\iota_1(x)\bigr)=-\iota_2(\bar x)$, and
$\tau\bigl(\iota_2(x)\bigr)=-\iota_1(\bar x)$.
\end{itemize}
\end{theorem}

\medskip

Given an algebra with involution $(A,\cdot,\bar\ )$, the
\emph{Steinberg unitary Lie algebra} $\frstu_3(A,\cdot,\bar\ )$ is
defined as the Lie algebra generated by the symbols $u_{ij}(a)$,
$1\leq i\ne j\leq 3$, $a\in A$, subject to the relations:
\[
\begin{split}
&u_{ij}(a)=u_{ji}(-\bar a),\\
&a\mapsto u_{ij}(a)\ \text{is linear,}\\
&[u_{ij}(a),u_{jk}(b)]=u_{ik}(ab)\ \text{for distinct $i,j,k$.}
\end{split}
\]
 Then \cite[Lemma 1.1]{AF93},
\[
\frstu_3(A,\cdot,\bar\ )=\frs\oplus u_{12}(A)\oplus u_{23}(A)\oplus
u_{31}(A),
\]
with $\frs=\sum_{i<j}[u_{ij}(A),u_{ij}(A)]$. This is a
$\bZ_2\times\bZ_2$-grading of $\frstu_3(A,\cdot,\bar\ )$.

\begin{proposition}\label{pr:steinberg}
Let $(A,\cdot,\bar\ )$ be a nonzero normal LRTA, then there is a
surjective homomorphism of Lie algebras
\[
\phi: \frstu_3(A,\cdot,\bar\ )\longrightarrow \frg(A,\cdot,\bar\ )
\]
such that, for any $a\in A$,
$\phi\bigl(u_{12}(a)\bigr)=-\iota_0(a)$,
$\phi\bigl(u_{23}(a)\bigr)=-\iota_1(a)$, and
$\phi\bigl(u_{31}(a)\bigr)=-\iota_2(a)$.
\end{proposition}
\begin{proof}
It is enough to realize that
$\bigl[u_{12}(a),u_{23}(b)\bigr]=u_{13}(a\cdot
b)=-u_{31}(\overline{a\cdot b})$, while
$\bigl[\iota_0(a),\iota_1(b)\bigr]=\iota_2(\overline{a\cdot b})$,
and cyclically.
\end{proof}

Also, the following result is proved by a straightforward
computation:

\begin{proposition}\label{pr:S4steinberg}
Let $(A,\cdot,\bar\ )$ be an algebra with involution. Then
$\frstu_3(A,\cdot,\bar\ )$ is endowed with an action of $S_4$ by
automorphisms by means of:
\begin{align*}
\tau_1 : u_{12}(a)&\mapsto u_{12}(a)& \tau_2 : u_{12}(a)&\mapsto
-u_{12}(a)\\
u_{23}(a)&\mapsto -u_{23}(a)&u_{23}(a)&\mapsto u_{23}(a)\\
u_{31}(a)&\mapsto -u_{31}(a)&u_{31}(a)&\mapsto -u_{31}(a)\\[10pt]
\varphi : u_{12}(a)&\mapsto u_{23}(a)& \tau : u_{12}(a)&\mapsto -u_{12}(\bar a)\\
u_{23}(a)&\mapsto u_{31}(a)&u_{23}(a)&\mapsto -u_{31}(\bar a)\\
u_{32}(a)&\mapsto u_{12}(a)&u_{31}(a)&\mapsto -u_{23}(\bar a)
\end{align*}
\end{proposition}

With these last two results, the result in \cite[Theorem 2.6]{Oku05}
follows easily:

\begin{theorem}\label{th:unitalstructurable}
The unital normal LRTA's are precisely the structurable algebras.
\end{theorem}
\begin{proof}
By Proposition \ref{pr:steinberg}, any unital normal LRTA is
$3$-faithful, and hence it is a structurable algebra \cite[Theorem
5.5]{AF93}. Conversely, if $(A,\cdot,\bar\ )$ is a structurable
algebra, then $\frg=\frstu_3(A,\cdot,\bar\ )$ is endowed with an
action of $S_4$ by automorphisms (Proposition \ref{pr:S4steinberg}),
which shows that $(A,\cdot,\bar\ )$ is a normal LRTA because of
Theorem \ref{th:deltasNLRTA}.
\end{proof}

\medskip

With the same arguments as for Theorem \ref{th:simpleSTA} one gets
too:

\begin{theorem}
Let $(A,\cdot,\bar\ )$ be a normal Lie related triple algebra with
$A\cdot A\ne 0$. Then $(A,\cdot,\bar\ )$ is simple if and only if
$\frg(A,\cdot,\bar\ )$ is simple as a Lie algebra with $S_4$-action.
\end{theorem}

\bigskip

\section{Examples}

\begin{example}\label{ex:allison}\quad (Structurable algebras)

Theorem \ref{th:unitalstructurable} shows that structurable algebras
are precisely the unital normal LRTA's, thus providing many examples
of these latter algebras. Given any structurable algebra
$(A,\cdot,\bar\ )$, the associated Lie algebra with $S_4$-action in
Theorem \ref{th:reverseLRTA} is
\[
\frg(A,\cdot,\bar\ )=\frt\oplus \bigl(\oplus_{i=0}^2
\iota_i(A)\bigr),
\]
where $\frt=\sum_{i=0}^2\theta^i\bigl(\delta(A,A)\bigr)$. But
Theorem \ref{th:deltasNLRTA} shows that
\[
\begin{split}
\delta_0(x,y)&=\delta_0(x,y\cdot 1)=
  -\delta_1(\bar y,\bar x)-\delta_2(1,\bar x\cdot y)\\
  &=-\bigl(L_xL{\bar y}-L_yL_{\bar x}\bigr)
   -\bigl(R_{\bar y\cdot x}-R_{\bar x\cdot y}\bigr)\\
  &=R_{(\bar x\cdot y-\bar y\cdot x)}+L_yL_{\bar x}-L_xL_{\bar y}.
\end{split}
\]
Therefore, $\frt$ is precisely the subspace $\calT_I$ of \emph{inner
triples} in \cite[Equation~(I)]{AF93}.

By identifying $u_{12}(x)$ with $-\iota_0(x)$, $u_{23}(x)$ with
$-\iota_1(x)$ and $u_{31}(x)$ with $-\iota_2(x)$ as in Proposition
\ref{pr:steinberg}, it follows that $\frg(A,\cdot,\bar\ )$ is
precisely the Lie algebra $\calK\bigl(A,\bar\ ,\gamma,\frt\bigr)$
constructed in \cite[Section~4]{AF93}, with $\gamma=(1,1,1)$. (Note
that other choices of $\gamma$ prevent this algebra from having the
symmetry induced by the action of $S_4$.) \qed
\end{example}

\bigskip

\begin{example}\label{ex:Jordan}\quad (Jordan algebras)

Unital Jordan algebras are examples of structurable algebras (where
the involution is the identity map). Nonunital Jordan algebras are
examples too of normal LRTA's.

Given a Jordan algebra $J$, with multiplication $\cdot$ and
involution $\bar\ =id$, the left and right multiplications coincide:
$L_x=R_x$ for any $x$, and $J$ is a normal LRTA \cite[Example
2.3]{Oku05} with
\[
\delta_i(x,y)=-[L_x,L_y]
\]
for any $x,y\in J$ and $i\in\bZ_3$. Since $[L_x,L_y]$ is a
derivation of $J$, conditions (i)--(vi) in Theorem \ref{th:deltas}
are clearly satisfied. Hence,
\[
\frt=\sum_{i=0}^2\theta^i\delta(J,J)
=\espan{\bigl([L_x,L_y],[L_x,L_y],[L_x,L_y]\bigr): x,y\in J},
\]
which is isomorphic to $\inder(J)=\espan{[L_x,L_y]: x,y\in J}$, the
Lie algebra of inner derivations (see \cite{JacobsonJordan}). Thus
$\frg(J,\cdot,\bar\ )$ is isomorphic to the Lie algebra
\begin{equation}\label{eq:tits}
\frg=\inder(J)\oplus\bigl(\frs\otimes_FJ\bigr)
\end{equation}
considered in \cite{Tit62}, where $\frs$ is the three-dimensional
simple Lie algebra with a basis $\{e_0,e_1,e_2\}$ such that
$[e_i,e_{i+1}]=e_{i+2}$, indices modulo $3$. (Over the reals this is
just $\frsu_2$.)

The Lie bracket in $\frg$ is given, for any $d\in\inder(J)$,
$s,t\in\frs$ and $x,y\in J$, by
\[
\left\{\begin{aligned}
 &[d,s\otimes x]=s\otimes d(x),\\
 &[s\otimes x,t\otimes y]=[s,t]\otimes x\cdot
 y-\tfrac{1}{2}\kappa(s,t)[L_x,L_y],
 \end{aligned}\right.
\]
where $\kappa$ denotes the Killing form of $\frs$. The isomorphism
just sends $\iota_i(x)$ to $e_i\otimes x$ for any $x$ and $i$.

\smallskip

If $J$ is any unital Jordan algebra, its Tits-Kantor-Koecher Lie
algebra $\mathcal{TKK}(J)$ (see \cite[Ch.~VIII]{JacobsonJordan}) is
isomorphic to the Lie algebra in \eqref{eq:tits}, but with $\frs$
replaced by $\frsl_2$. If $\sqrt{-1}\in F$, then $\frs$ is
isomorphic to $\frsl_2$ and, hence, $\mathcal{TKK}(J)$ becomes a Lie
algebra with $S_4$-action. \qed
\end{example}

\bigskip

\begin{example}\quad (Lie algebras)

Given any Lie algebra $\frg$, the direct sum of four copies of
$\frg$: $\frg^4$, is endowed with a natural action of $S_4$:
\[
\begin{split}
S_4&\longrightarrow \Aut(\frg^4)\\
\sigma&\mapsto \Bigl( (x_1,x_2,x_3,x_4)\mapsto
 (x_{\sigma^{-1}(1)},x_{\sigma^{-1}(2)},x_{\sigma^{-1}(3)},
   x_{\sigma^{-1}(4)})\Bigr).
\end{split}
\]
The Lie algebra $\frg^4$ may be identified naturally with
$\frg\otimes_FF^4$, and $F^4$ contains the basis
$\{1,e_0,e_1,e_2\}$, where $1=(1,1,1,1)$, $e_0=(1,1,-1,-1)$,
$e_1=(1,-1,-1,1)$ and $e_2=(1,-1,1,-1)$. Then, the components in
\eqref{eq:tg0g1g2} become:
\[
\begin{split}
\frt&=\{X\in\frg^4:\tau_1(X)=X,\,\tau_2(X)=X\}=\frg\otimes 1,\\
\frg_0&=\{X\in\frg^4:\tau_1(X)=X,\,\tau_2(X)=-X\}=\frg\otimes e_0,\\
\frg_1&=\{X\in\frg^4:\tau_1(X)=-X,\,\tau_2(X)=X\}=\frg\otimes e_1,\\
\frg_2&=\{X\in\frg^4:\tau_1(X)=-X,\,\tau_2(X)=-X\}=\frg\otimes e_2.
\end{split}
\]
It follows that the corresponding normal LRTA can be identified with
$\frg$, with $x\cdot y=[x,y]$, $\bar x=-x$ and
$\delta(x,y)=\bigl(\ad_{[x,y]},\ad_{[x,y]},\ad_{[x,y]}\bigr)$ for
any $x,y\in \frg$, in accordance to \cite[Example 2.3]{Oku05}.
Therefore, any Lie algebra is a normal LRTA. \qed
\end{example}

\bigskip

\begin{example}\label{ex:LTS}\quad (Lie triple systems)

Any Lie triple system is the odd part of a $\bZ_2$-graded Lie
algebra $\frg=\frg\subo\oplus\frg\subuno$, with the triple product
being given by $[xyz]=[[x,y],z]$ for any $x,y,z\in\frg\subuno$.

Given a $\bZ_2$-graded Lie algebra
$\frg=\frg\subo\oplus\frg\subuno$, let us denote by $\nu$ the
grading automorphism: $\nu(x\subo+x\subuno)=x\subo-x\subuno$. Then
$S_4$ acts as automorphisms of $\frg^3$ as follows:
\[
\left\{\begin{aligned}
 \tau_1(x,y,z)&=(x,\nu(y),\nu(z)),\\
 \tau_2(x,y,z)&=(\nu(x),y,\nu(z)),\\
 \varphi(x,y,z)&=(z,x,y),\\
 \tau(x,y,z)&=(x,z,y).
 \end{aligned}\right.
\]
The subspaces in \eqref{eq:tg0g1g2} are
\[
\begin{split}
\frt&=\frg\subo^{\,3},\\
\frg_0&=\{(x\subuno,0,0): x\subuno\in\frg\subuno\},\\
\frg_1&=\{(0,x\subuno,0): x\subuno\in\frg\subuno\},\\
\frg_2&=\{(0,0,x\subuno): x\subuno\in\frg\subuno\}.
\end{split}
\]
Hence the associated normal LRTA is $A\simeq \frg\subuno$, with
$x\cdot y=0$, $\bar x=-x$ and
$\delta(x,y)=\bigl(\ad_{[x,y]},0,0\bigr)$ for any
$x,y\in\frg\subuno$.

If $\frg=\frg\subo\oplus\frg\subuno$ is a simple Lie algebra, then
$[\frg\subuno,\frg\subuno]=\frg\subo$ and $x\subo\mapsto
\ad_{x\subo}\vert_{\frg\subuno}$ is one-to-one. Hence, with
$A=\frg\subuno$ as before, $\frg(A,\cdot,\bar\ )$ coincides with
$\frg^3$, and this has no proper ideals invariant under the action
of $S_4$. Hence $\frg(A,\cdot,\bar\ )$ is simple as a Lie algebra
with $S_4$-action, even though $(A,\cdot)$ is a trivial algebra.
\qed
\end{example}

\bigskip

\begin{example}\quad (Tensor products of symmetric composition
algebras)

A \emph{symmetric composition algebra} is a triple $(S,*,q)$,  where
$q:S\rightarrow F$ is a regular quadratic form satisfying for any
$x,y,z\in S$:
\[
\begin{split}
 q(x*y)&=q(x)q(y),\\
 q(x*y,z)&=q(x,y*z),
\end{split}
\]
where $q(x,y)=q(x+y)-q(x)-q(y)$ is the polar of $q$ (see
\cite[Ch.~VIII]{KMRT98})

The classification of the symmetric composition algebras was
obtained in \cite{EM93}, \cite{Eld97}.

Given any unital composition algebra (or Hurwitz algebra) $C$ with
norm $q$ and standard involution $x\mapsto \bar x$, the new algebra
defined on $C$ but with multiplication
$$
x*y=\bar x\bar y,
$$
is a symmetric composition algebra, called the associated
\emph{para-Hurwitz algebra}. In dimensions $1$, $2$ or $4$, any
symmetric composition algebra is a para-Hurwitz algebra (with a few
exceptions in dimension $2$), while in dimension $8$, apart from the
para-Hurwitz algebras, there is a new family of symmetric
composition algebras termed \emph{Okubo algebras}.

If $(S,*,q)$ is any symmetric composition algebra, consider the
corresponding orthogonal Lie algebra $\fro(S,q)=\{ d\in \frgl(S):
q\bigl(d(x),y\bigr)+q\bigl(x,d(y)\bigr)=0\ \forall x,y\in S\}$, and
the subalgebra of $\fro(S,q)^3$ defined by
\[
\begin{split}
 \tri&(S,*,q)\\
  &=\{(d_0,d_1,d_2)\in \fro(S,q)^3 :
 d_0(x*y)=d_1(x)*y+x*d_2(y)\ \forall x,y\in S\}\\
 &=\{(d_0,d_1,d_2)\in \fro(S,q)^3 :
 \langle d_0(x),y,z\rangle+\langle x,d_1(y),z\rangle
 +\langle x,y,d_2(z)\rangle=0\\
 &\hspace{3.5in} \forall x,y,z\in S\},
\end{split}
\]
where $\langle x,y,z\rangle =q(x,y*z)$. It turns out that the map,
\[
\begin{split}
 \theta: \tri(S,*,q)&\rightarrow \tri(S,*,q)\\
 (d_0,d_1,d_2)&\mapsto (d_2,d_0,d_1)
\end{split}
\]
is an automorphism of $(S,*,q)$ of order $3$. Its fixed subalgebra
is (isomorphic to) the derivation algebra of $(S,*)$ which, if the
dimension is $8$ and the characteristic of the ground field is $\ne
2,3$, is a simple Lie algebra of type $G_2$ in the para-Hurwitz case
and a simple Lie algebra of type $A_2$ (a form of $\frsl_3$) in the
Okubo case.

\smallskip

A straightforward computation (see \cite{EO01} for a more general
setting) shows that for any $x,y\in S$, the triple
\[
t_{x,y}=
 \left( \sigma_{x,y},\frac{1}{2}q(x,y)id-r_xl_y,
 \frac{1}{2}q(x,y)id-l_xr_y\right)
\]
is in $\tri(S,*,q)$, where $\sigma_{x,y}(z)=q(x,z)y-q(y,z)x$,
$r_x(z)=z*x$, and $l_x(z)=x*z$ for any $x,y,z\in S$.

\medskip

Given two symmetric composition algebras $(S,*,q)$ and $(S',*,q')$,
their tensor product $A=S\otimes_FS'$ is a normal STA with
multiplication
\[
(a\otimes x)*(b\otimes y)=(a*b)\otimes (x*y),
\]
 and skew-symmetric bilinear map
$\delta:A\times A\rightarrow \stri(A,*)$ given by
\[
\delta(a\otimes x,b\otimes y)
 =q'(x,y)t_{a,b}+q(a,b)t'_{x,y}\in\tri(S,*,q)\oplus\tri(S',*,q')
\subseteq \stri(A,*),
\]
for any $a,b\in S$ and $x,y\in S'$ (see \cite{Eld04}).

\smallskip

This example can be extended by considering the so called
\emph{generalized symmetric composition algebras} \cite[Examples 2.4
and 2.5]{Oku05}.

Note also that the tensor product of two unital composition algebras
is a structurable algebra, and hence a normal LRTA. \qed
\end{example}

\bigskip

\begin{example}

Any Okubo algebra $(S,*,q)$ over $F$ either contains a nonzero
idempotent, or there exists a cubic field extension $L/F$ such that
$L\otimes_FS$ contains a nonzero idempotent. In the first case,
there is an order $3$ automorphism $\phi$  of $S$ such that $S$
becomes a para-Hurwitz algebra with the new multiplication given by
$x\star y=\phi(x)*\phi^2(y)$, for any $x,y$, and the same norm $q$
(see \cite[\S 3.4]{KMRT98}). This is used in \cite{Eld05'} to relate
the models of the exceptional Lie algebras in \cite{Eld04}
constructed in terms of Okubo algebras, to those given by
para-Hurwitz algebras.

This can be extended as follows. Let $(A,*)$ be a normal STA and
consider an order $3$ automorphism $\phi\in\Aut(A,*)$  such that
$\phi\delta_0(x,y)\phi^{-1}=\delta_0\bigl(\phi(x),\phi(y)\bigr)$ for
any $x,y\in A$.

Note that since $l_{\phi(x)}=\phi l_x\phi^{-1}$ and
$r_{\phi(x)}=\phi r_x\phi^{-1}$, because $\phi$ is an automorphism,
then
$\phi\delta_i(x,y)\phi^{-1}=\delta_i\bigl(\phi(x),\phi(y)\bigr)$ for
$i=1,2$, because of Theorem \ref{th:deltas}. Also, if
$(d_0,d_1,d_2)\in\stri(A,*)$, then $\bigl(\phi d_0\phi^{-1},\phi
d_1\phi^{-1},\phi d_2\phi^{-1}\bigr)\in\stri(A,*)$ too, and hence,
for any $x,y\in A$:
\[
\left\{\begin{aligned}
 &\Bigl(\delta_0\bigl(\phi(x),\phi(y)\bigr),
   \delta_1\bigl(\phi(x),\phi(y)\bigr),\delta_1\bigl(\phi(x),\phi(y)\bigr)\Bigr)
   \in\stri(A,*),\\
 &\Bigl(\phi\delta_0(x,y)\phi^{-1},\phi\delta_1(x,y)\phi^{-1},
   \phi\delta_2(x,y)\phi^{-1}\Bigr)\in\stri(A,*),
 \end{aligned}\right.
\]
so its difference also belongs to $\stri(A,*)$:
\[
\Bigl(\phi\delta_0(x,y)\phi^{-1} -
   \delta_0\bigl(\phi(x),\phi(y)\bigr),\, 0\, ,\, 0\,
   \Bigr)\in\stri(A,*).
\]
But if either $A=A*A$, or $\{x\in A: x*A=0\}=0$, or $\{x\in A:
A*x=0\}=0$, then $(d_0,0,0)\in\stri(A,*)$ implies that $d_0=0$.
Hence the condition $\phi\delta_0(x,y)\phi^{-1} =
   \delta_0\bigl(\phi(x),\phi(y)\bigr)$ is superfluous in these
cases.

Define a new multiplication on $A$ by means of
\[
x\star y=\phi(x)*\phi^2(y)
\]
for any $x,y\in A$ as before. Then it is easily checked that
$(A,\star)$ is again a normal STA with $\delta^\star:A\times
A\rightarrow \stri(A,\star)$ given by
\[
\delta_i^\star (x,y)=\phi^{-i}\delta_i(x,y)\phi^i
\]
for any $x,y\in A$ and $i\in\bZ_3$, because
\[
\begin{split}
\delta_i^\star(x,y)(u\star v)
 &=\phi^{-i}\delta_i(x,y)\phi^i\bigl(\phi(u)*\phi^2(v)\bigr)\\
 &=\phi^{-i}\delta_i(x,y)\bigl(\phi^{i+1}(u)*\phi^{i+2}(v)\bigr)\\
 &=\phi^{-i}\Bigl(\delta_{i+1}(x,y)\bigl(\phi^{i+1}(u)\bigr)*\phi^{i+2}(v)
   \\
   &\qquad\qquad\qquad+\phi^{i+1}(u)*\delta_{i+2}(x,y)\bigl(\phi^{i+2}(v)\bigr)\Bigr)\\
 &=\phi\bigl(\delta_{i+2}^\star(x,y)(u)\bigr)*\phi^2(v)+
   \phi(u)*\phi^2\bigl(\delta_{i+2}^\star(x,y)(v)\bigr)\\
 &=\delta_{i+2}^\star(x,y)(u)\star v+
   u \star \delta_{i+2}^\star(x,y)(v),
\end{split}
\]
for any $x,y,u,v\in A$ and $i\in\bZ_3$.

Under these conditions, the linear map:
\[
\Phi: \frg(A,\star)\longrightarrow \frg(A,*),
\]
such that
\[
\left\{\begin{aligned}
 &\Phi\bigl(\iota_i^\star(x)\bigr)=\iota_i\bigl(\phi^i(x)\bigr),\\
 &\Phi(d_0,d_1,d_2)= (d_0,\phi d_1\phi^{2},\phi^2d_2\phi),
 \end{aligned}\right.
\]
is an isomorphism of Lie algebras.

\smallskip

If now $(A,\cdot,\bar\ )$ is a normal LRTA and $\phi$ is an order
$3$ automorphism of $(A,\cdot)$ such that $\bar\phi=\phi^2$, then a
new multiplication can be defined by
\[
x\bullet y=\phi(x)\cdot\phi^2(y),
\]
for any $x,y\in A$, which satisfies
\[
\begin{split}
\overline{x\bullet y}&=\overline{\phi(x)\cdot \phi^2(y)}
 =\overline{\phi^2(x)}\cdot \overline{\phi(x)}\\
 &=\overline{\phi^2}(\bar y)\cdot\bar\phi (\bar x)=
  \phi(\bar y)\cdot\phi^2(\bar x)\\
  &=\bar y\bullet \bar x.
\end{split}
\]
Hence $\bar\ $ is an involution too of $(A,\bullet)$. Now, assuming
that
$\phi\delta_0(x,y)\phi^{-1}=\delta_0\bigl(\phi(x),\phi(y)\bigr)$ for
any $x,y\in A$, $(A,\bullet,\bar\ )$ is another normal LRTA relative
to the skew-symmetric bilinear map $\delta^\bullet:A\times
A\rightarrow \lrt(A,\bullet,\bar\ )$ given by:
\[
\delta_i^\bullet(x,y)=\phi^{-i}\delta_i(x,y)\phi^i,
\]
for any $x,y\in A$ and $i\in\bZ_3$.

\smallskip

This allows the construction of new normal STA's or LRTA's out of
old ones. For instance, let $\frs$ be the three-dimensional simple
Lie algebra that appears in Example \ref{ex:Jordan}. This is a
normal LRTA, since it is a Lie algebra. Consider the order $3$
automorphism given by $\phi(e_i)=e_{i+1}$, for any $i\in\bZ_3$. This
automorphism induces the new multiplication $\bullet$ given by
\[
s\bullet t=[\phi(s),\phi^2(t)]
\]
which gives:
\[
e_i\bullet e_{i+1}=-e_{i+2},\quad
 e_{i+1}\bullet e_i=0,\quad
 e_i\bullet e_i=e_i,
\]
for any $i\in\bZ_3$. This normal LRTA does not fit into any of the
previous examples. \qed
\end{example}

\bigskip

\section{$S_4$-action on Kantor's construction}

Given any structurable algebra $(A,\bar\ )$ (multiplication denoted
by juxtaposition), Kantor's construction gives a $5$-graded Lie
algebra $\calK(A,\bar\ ,\frd)$ (see \cite{All79}, based on
\cite{Kan72}, \cite{Kan73}), where $\frd$ is a Lie subalgebra of the
Lie algebra of derivations which contains the inner derivations. On
the other hand, Allison and Faulkner  defined a
$\bZ_2\times\bZ_2$-graded Lie algebra $\calK(A,\bar\ ,\gamma,\frv)$
in \cite[Section 4]{AF93}, where $\gamma\in F^3$, and $\frv$ is a
Lie subalgebra of $\lrt(A,\bar\ )$ containing the subspace of inner
triples (see Example \ref{ex:allison}). A precedent of this
construction appears in \cite{Vin66}.

The aim of this section is to show that the Lie algebra
$\calK(A,\bar\ ,\frd)$ is isomorphic to the Lie algebra
$\calK(A,\bar\ ,\gamma,\frv)$, for some suitable $\gamma$ and
$\frv$, and that, under some restrictions on the ground field, it is
endowed with an action of the symmetric group $S_4$.

\medskip

Let $(A,\bar\ )$ be a structurable algebra over a ground field of
characteristic $\ne 2,3$. This assumption on the field will be kept
throughout the section. Let $\inder(A,\bar\ )$ (the \emph{inner
derivation algebra}) be the linear span of the derivations $\{
D_{x,y}:x,y\in A\}$, where \cite[Eq.~(15)]{All78}:
\[
D_{x,y}(z)=\frac{1}{3}\bigl[[x,y]+[\bar x,\bar y],z\bigr]
 +(z,y,x)-(z,\bar x,\bar y).
\]
(Here $(x,y,z)=(xy)z-x(yz)$ is the associator.) The subspace
$\inder(A,\bar\ )$ is an ideal of the Lie algebra of derivations
$\der(A,\bar\ )$. Also, consider the Lie algebra of inner related
triples
\[
\inlrt(A,\bar\ )=\sum_{i=0}^2\espan{\theta^i\bigl(\delta(x,y)\bigr):
x,y\in A},
\]
where
\[
\delta(x,y) =
 \Bigl(R_{(\bar x y-\bar yx)}+L_yL_{\bar x}-L_xL_{\bar y},
 L_{\bar y}L_x-L_{\bar x}L_y, R_{\bar y}R_x-R_{\bar x}R_y\Bigr).
\]
(See Example \ref{ex:allison} and  note that $\inlrt(A,\bar\ )$ was
denoted by $\calT_I$ in \cite{AF93}.)

The following notation will be used. Given the structurable algebra
$(A,\bar\ )$, let $S$ be the set of skew-symmetric elements:
$S=\{x\in A: \bar x=-x\}$. Then $\calT_S$ will denote the vector
space
\[
\calT_S=\{
\bigl(L_{s_1}-R_{s_2},L_{s_2}-R_{s_0},L_{s_0}-R_{s_1}\bigr) :
s_0,s_1,s_2\in S,\ s_0+s_1+s_2=0\}.
\]
Also, for any subspace $\frs\subseteq \frgl(A)$, denote by
$\frs^{<3>}$ the `diagonal subspace' $\{(s,s,s): s\in\frs\}\leq
\frgl(A)^3$.

\begin{lemma}\label{le:lrt}
Let $(A,\bar\ )$ be a structurable algebra. Then:
\begin{alphaenumerate}
\item For any $x,y\in A$, $\sum_{i=0}^2\delta_i(x,y)=-3D_{\bar
x,y}$.
\item $\lrt(A,\bar\ )= \der(A,\bar\ )^{<3>}\oplus\calT_S$, while
$\inlrt(A,\bar\ )=\inder(A,\bar\ )^{<3>}\oplus\calT_S$.
\item The map
\[
\begin{split}
 \bigl\{ \frd: \inder(A,\bar\ )\leq \frd\leq \der(A,\bar\
)\bigr\} &\rightarrow \bigl\{ \frl: \inlrt(A,\bar\ )\leq \frl\leq
\lrt(A,\bar\ )\bigr\}\\
\frd\qquad&\mapsto\qquad\ \frd^{<3>}\oplus \calT_S,
\end{split}
\]
is a bijection.
\end{alphaenumerate}
\end{lemma}
\begin{proof} By \cite[Eq.~(A1)]{AF93},
\[
(x,\bar y,z)-(y,\bar x,z)=(z,\bar x,y)-(z,\bar y,x),
\]
which is equivalent to
\[
L_{(x\bar y-y\bar x)}+R_{(\bar xy-\bar yx)}=
 L_xL_{\bar y}-L_yL_{\bar x}-R_xR_{\bar y}+R_yR_{\bar x}.
\]
Hence
\[
\begin{split}
L_{([\bar x,y]+[x,\bar y])}&=
  L_{(x\bar y-y\bar x)}+L_{(\bar xy-\bar yx)}\\
  &=-R_{(\bar xy-\bar yx)}+\bigl(L_xL_{\bar y}-L_yL_{\bar x}\bigr))-
   \bigl(R_xR_{\bar y}-R_yR_{\bar x}\bigr)\\
   &\quad -R_{(x\bar y-y\bar x)} +
  \bigl(L_{\bar x}L_y-L_{\bar y}L_x\bigr)-
  \bigl(R_{\bar x}R_y-R_{\bar y}R_x\bigr).
\end{split}
\]
Thus,
\[
\begin{split}
3D_{\bar x,y}&=L_{([\bar x,y]+[x,\bar y])}
  -R_{([\bar x,y]+[x,\bar y])}
   +3R_{(x\bar y-y\bar x)}
   +3\bigl(R_{\bar x}R_y-R_{\bar y}R_x\bigr)\\
  &=-2\bigl(R_{(\bar xy-\bar yx)}+R_{(x\bar y-y\bar x)}\bigr)
     +3R_{(x\bar y-y\bar x)}
     +3\bigl(R_{\bar x}R_y-R_{\bar y}R_x\bigr)\\
  &\quad +
    \bigl(L_xL_{\bar y}-L_yL_{\bar x}\bigr) +
    \bigl(L_{\bar x}L_y-L_{\bar y}L_x\bigr)\\
  &\quad -\bigl(R_xR_{\bar y}-R_yR_{\bar x}\bigr)-
     \bigl(R_{\bar x}R_y-R_{\bar y}R_x\bigr)\\
  &=\bigl(R_{(x\bar y-y\bar x)}-2R_{(\bar xy-\bar yx)}\bigr)
    -\bigl(R_xR_{\bar y}-R_yR_{\bar x})+2\bigl(R_{\bar x}R_y-R_{\bar
    y}R_x\bigr)\\
  &\quad +\bigl(L_xL_{\bar y}-L_yL_{\bar x}\bigr)
     +\bigl(L_{\bar x}L_y-L_{\bar y}L_x\bigr).
\end{split}
\]
But $D_{\bar x,y}=D_{x,\bar y}$ \cite[Lemma 6]{All78}, so
\[
\begin{split}
9D_{\bar x,y}&=6D_{\bar x,y}+3D_{x,\bar y}\\
 &= -3R_{(\bar xy-\bar yx)}+3\bigl(R_{\bar x}R_y-R_{\bar
 y}R_x\bigr)\\
 &\quad +3\bigl(L_xL_{\bar y}-L_yL_{\bar x}\bigr)
   +3\bigl(L_{\bar x}L_y-L_{\bar y}L_x\bigr),
\end{split}
\]
which proves (a).

Now, (b) follows from \cite[Corollary 3.5]{AF93} and its proof,
using (a).

For (c) consider the linear map
\[
\begin{split}
\phi: \lrt(A,\bar\ )&\longrightarrow S\times S\\
(d_0,d_1,d_2)&\mapsto \bigl(d_1(1),d_2(1)\bigr).
\end{split}
\]
Any $(d_0,d_1,d_2)\in\ker\phi$ satisfies $\bar d_0(1)=0$, so
$d_0(1)=0$, and
\[
\bar d_i(x)= \begin{cases} \bar d_i(x1)=d_{i+1}(x)\\ \bar
d_i(1x)=d_{i+2}(x)\end{cases}
\]
for any $x\in A$, so $\bar d_i=d_{i+1}=d_{i+2}$ for any $i$. This
implies that $d_i=d_{i+1}=d_{i+2}=\bar d_i$ for any $i$, and hence
$\ker\phi=\der(A,\bar\ )^{<3>}$. It follows that $\lrt(A,\bar\
)=\ker\phi\oplus\calT_S$, $\inlrt(A,\bar\
)=\ker\phi\vert_{\inlrt(A,\bar\ )}\oplus\calT_S$, and for any
subalgebra $\frl$ of $lrt(A,\bar\ )$ containing $\inlrt(A,\bar\ )$,
$\frl=\ker\phi\vert_{\frl}\oplus\calT_S$. Finally, the map
$\frl\mapsto \ker\phi\vert_{\frl}$ is the inverse of the map in (c).
\end{proof}

\smallskip

For any subalgebra $\frd$ with $\inder(A,\bar\ )\leq \frd\leq
\der(A,\bar\ )$, denote by $\frl(A,\bar\ ,\frd)$ the subalgebra
$\frd^{<3>}\oplus\calT_S$ of $\lrt(A,\bar\ )$.

\medskip

Let us recall Kantor's construction from \cite[Section 3]{All79}.

Let $\frd$ be a subalgebra of $\der(A,\bar\ )$ containing
$\inder(A,\bar\ )$, then
\[
\calK(A,\bar\ ,\frd)=\tilde \frn\oplus
\bigl(T_A\oplus\frd\bigr)\oplus \frn,
\]
where $\frn=A\times S$, $\tilde\frn$ is another copy of $\frn$, and
$T_x=V_{x,1}$, where $V_{x,y}(z)=(x\bar y)z+(z\bar y)x-(z\bar x)y$
for any $x,y,z\in A$. Then $T_A\oplus\frd$ is a Lie subalgebra of
$\frgl(A)$ and the bracket of any two of its elements in
$\calK(A,\bar\ ,\frd)$ coincides with its bracket in $\frgl(A)$.
Moreover,
\[
\left\{\begin{aligned} &[f,(x,s)]=\bigl(f(x),f^{\delta}(s)\bigr),\\
   &[f,(x,s)\tilde\
   ]=\bigl(f^{\epsilon}(x),f^{\epsilon\delta}(x)\bigr)\tilde\ ,\\
   &[(x,r),(y,s)]=(0,x\bar y-y\bar x),\\
   &[(x,r)\tilde\ ,(y,s)\tilde\ ]=(0,x\bar y -y\bar x)\tilde\ ,\\
   &[(x,r),(y,s)\tilde\ ]=-(sx,0)\tilde\ +V_{x,y}+L_rL_s+(ry,0),
\end{aligned}\right.
\]
for any $x,y\in A$, $r,s\in S$, $f\in T_A\oplus\frd$, where
$f^\epsilon=f-T_{(f(1)+\overline{f(1)})}$ and
$f^\delta=f+R_{\overline{f(1)}}$, so that $d^\epsilon=d^\delta=d$
for any $d\in\der(A,-)$, while $T_x^\epsilon=-T_{\bar x}$ and
$T_x^\delta=T_x+R_{\overline{2x-\bar x}}=L_x+R_{\bar x}$ for any
$x\in A$.

\smallskip

The Lie algebra $\calK =\calK(A,\bar\ ,\frd)$ is $5$-graded:
\[
\calK=\calK_{-2}\oplus\calK_{-1}\oplus\calK_0\oplus\calK_1\oplus\calK_2,
\]
where $\calK_{-2}=(0\times S)\tilde\ $, $\calK_{-1}=(A\times
0)\tilde\ $, $\calK_0=T_A\oplus\frd$, $\calK_1=A\times 0$ and
$\calK_2=0\times S$. Allison \cite[Eq.~(15)]{All79} considered the
order two automorphism $\chi$ such that
\[
\chi(\tilde m+f+n)=\tilde n+f^\epsilon+m,
\]
for any $m,n\in A\times
S$ and $f\in T_A\oplus\frd$. Also, for any $0\ne\alpha\in F$,
$\sigma_\alpha(x_j)=\alpha^jx_j$ ($x_j\in\calK_j$) defines an
automorphism too.

Fix $0\ne\alpha\in F$ and consider the two commuting order two
automorphisms $\tau_1,\tau_2$ of $\calK(A,\bar\ ,\frd)$ given by
\[
\tau_1=\sigma_{-1},\qquad \tau_2=\sigma_{\alpha}\chi.
\]
These two automorphisms induce a $\bZ_2\times\bZ_2$-grading of
$\calK=\calK(A,\bar\ ,\frd)$ as in \eqref{eq:tg0g1g2}, where
\[
\begin{split}
\calK_{(\bar 0,\bar 0)}&=T_S\oplus\frd
     \oplus\{\alpha(0,s)+\alpha^{-1}(0,s)\tilde\ : s\in S\},\\
\calK_{(\bar 1,\bar 0)}&=T_H
      \oplus\{\alpha(0,s)-\alpha^{-1}(0,s)\tilde\ : s\in S\},\\
\calK_{(\bar 0,\bar 1)}&=\{(x,0)+\alpha^{-1}(x,0)\tilde\ :x\in
A\},\\
\calK_{(\bar 1,\bar 1)}&=\{\alpha(x,0)-(x,0)\tilde\ :x\in A\}.
\end{split}
\]
Note that $\calK_{(\bar 1,\bar 0)}$, $\calK_{(\bar 0,\bar 1)}$ and
$\calK_{(\bar 1,\bar 1)}$ are vector spaces isomorphic to $A$.

For ease of notation, write
\[
\left\{\begin{aligned}
 \epsilon_1(x)&=(x,0)+\alpha^{-1}(x,0)\tilde\ \in\calK_{(\bar 0,\bar 1)},\\
 \epsilon_2(x)&=\alpha(\bar x,0)-(\bar x,0)\tilde\ ,\in\calK_{(\bar 1,\bar 0)}\\
 \epsilon_0(x)&=\frac{1}{2}\Bigl(T_{(x+\bar
 x)}+\bigl(\alpha(0,x-\bar x)-\alpha^{-1}(0,x-\bar x)\tilde\
 \bigr)\Bigr)\in\calK_{(\bar 1,\bar 1)}.
\end{aligned}\right.
\]
Then the bracket in $\calK$ gives, for any $x,y\in A$:
\[
\begin{split}
[\epsilon_1(x),\epsilon_2(y)]
 &= \bigl[(x,0)+\alpha^{-1}(x,0)\tilde\ ,
   \alpha(\bar y,0)-(\bar y,0)\tilde\ \bigr]\\
 &=\alpha(0,xy-\bar y\bar x)-\alpha^{-1}(0,xy-\bar y\bar x)\tilde\
    -V_{x,\bar y}-V_{\bar y,x}\\
 &=\alpha(0,xy-\overline{xy},0)-
      \alpha^{-1}(0,xy-\overline{xy})\tilde\
     -T_{(xy+\overline{xy})}\\
 &=-2\epsilon_0(\overline{xy}),
\end{split}
\]
as $\bigl(V_{x,y}+V_{y,x}\bigr)(z)=(x\bar y+y\bar x)z=T_{(x\bar
y+y\bar x)}(z)$.

Also, for any $a\in H\,(=\{ x\in A: \bar x=x\})$, $x\in A$ and $s\in
S$, $[T_a,(x,0)]=(ax,0)$, while $[T_a,(x,0)\tilde\ ]=-(ax,0)\tilde\
$, and
\[
[\alpha(0,s)-\alpha^{-1}(0,s)\tilde\
,(x,0)]=\alpha^{-1}[(x,0),(0,s)]=-\alpha^{-1}(sx,0)\tilde\ ,
\]
while
\[
[\alpha(0,s)-\alpha^{-1}(0,s)\tilde\ ,(x,0)\tilde\
]=\alpha[(0,s),(x,0)\tilde\ ]=\alpha(sx,0).
\]
Therefore,
\[
\begin{split}
[\epsilon_0(x),\epsilon_1(y)]&=
  \frac{1}{2}\Bigl[T_{(x+\bar x)}+\alpha(0,x-\bar
  x)-\alpha^{-1}(0,x-\bar x)\tilde\ , (y,0)+\alpha^{-1}(y,0)\tilde\
  \Bigr]\\
  &=\frac{1}{2}\Bigl( ((x+\bar x)y,0)-\alpha^{-1}((x+\bar
  x)y,0)\tilde\ \\
  &\qquad\qquad\qquad-\alpha^{-1}((x-\bar x)y,0)\tilde\ +((x-\bar
  x)y,0)\Bigr)\\
  &=(xy,0)-\alpha^{-1}(xy,0)\tilde\ \\
  &=\alpha^{-1}\epsilon_2(\overline{xy}),
\end{split}
\]
and
\[
\begin{split}
[\epsilon_2(x),\epsilon_0(y)]&=
  -\frac{1}{2}\Bigl[T_{(y+\bar y)}+\alpha(0,y-\bar
  y)-\alpha^{-1}(0,y-\bar y)\tilde\ , \alpha(\bar x,0)-(\bar x,0)\tilde\
  \Bigr]\\
  &=-\frac{1}{2}\Bigl( \alpha((y+\bar y)\bar x,0)+((y+\bar
  y)\bar x,0)\tilde\ \\
  &\qquad\qquad\qquad-((y-\bar y)\bar x,0)\tilde\ -\alpha((y-\bar
  y)\bar x,0)\Bigr)\\
  &=-\alpha(\bar y\bar x,0)+(\bar y\bar x,0)\tilde\ \\
  &=-\alpha\epsilon_1(\overline{xy}).
\end{split}
\]
That is,
\begin{equation}\label{eq:epsilones}
\begin{split}
[\epsilon_0(x),\epsilon_1(y)]&=\alpha^{-1}\epsilon_2(\overline{xy}),\\
[\epsilon_1(x),\epsilon_2(y)]&=-2\epsilon_0(\overline{xy}),\\
[\epsilon_2(x),\epsilon_0(y)]&=-\alpha\epsilon_1(\overline{xy}).
\end{split}
\end{equation}

\smallskip

\begin{lemma}\label{le:K00}
The linear map
\[
\begin{split}
\psi:\calK_{(\bar 0,\bar 0)}&\longrightarrow \lrt(A,\bar\ )\\
p\quad&\mapsto \bigl(\delta_0(p),\delta_1(p),\delta_2(p)\bigr),
\end{split}
\]
where $\delta_i(p)$ is determined by
\[
\bigl[p,\epsilon_i(x)\bigr]=\epsilon_i\bigl(\delta_i(p)(x)\bigr),
\]
is a one-to-one Lie algebra homomorphism with image $\frl(A,\bar\
,\frd)$.
\end{lemma}
\begin{proof}
First of all, $\psi$ is well defined since the Jacobi identity shows
that, for any $p\in\calK_{(\bar 0,\bar 0)}$, $x,y\in A$ and
$i\in\bZ_3$,
\[
\begin{split}
\bigl[p,\bigl[\epsilon_{i+1}(x)&,\epsilon_{i+2}(y)\bigr]\bigr]\\
 &=\bigl[\bigl[p,\epsilon_{i+1}(x)\bigr],\epsilon_{i+2}(y)\bigr]
  +\bigl[\epsilon_{i+1}(x),\bigl[p,\epsilon_{i+2}(y)\bigr]\bigr]\\
 &=\bigl[\epsilon_{i+1}\bigl(\delta_{i+1}(p)(x)\bigr),\epsilon_{i+2}(y)\bigr]
  +\bigl[\epsilon_{i+1}(x),\epsilon_{i+2}\bigl(\delta_{i+2}(p)(y)\bigr)\bigr]\\
 &=\xi_i\epsilon_i\bigl(\overline{\delta_{i+1}(p)(x)y+x\delta_{i+2}(p)(y)}\bigr),
\end{split}
\]
where $\xi_0=-2$, $\xi_1=-\alpha$, $\xi_2=\alpha^{-1}$ by
\eqref{eq:epsilones}, but also
\[
\bigl[p,\bigl[\epsilon_{i+1}(x),\epsilon_{i+2}(y)\bigr]\bigr]
 =\bigl[p,\xi_i\epsilon_i(\overline{xy})\bigr]
 =\xi_i\epsilon_i\bigl(\delta_i(p)(\overline{xy})\bigr),
\]
so
\[
\overline{\delta_i(p)}(xy)=\delta_{i+1}(p)(x)y+x\delta_{i+2}(p)(y).
\]
The fact that $\psi$ is a Lie algebra homomorphism is clear.

Besides, $\psi(\frd)=\frd^{<3>}$, and for any $s\in S$ and $x\in A$,
\[
\begin{split}
\bigl[T_s,\epsilon_1(x)\bigr]&=
 \bigl[T_s,(x,0)+\alpha^{-1}(x,0)\tilde\ \bigr]
 = (T_s(x),0)+\alpha^{-1}(T_s(x),0)\tilde\ \\
 &=\epsilon_1\bigl(T_s(x)\bigr),\\
\bigl[T_s,\epsilon_2(x)\bigr]&=
 \bigl[T_s,\alpha(\bar x,0)-(\bar x,0)\tilde\ \bigr]
 =\alpha(T_s(\bar x),0)-(T_s(\bar x),0)\tilde\ \\
 &=\alpha\bigl(\overline{\bar T_s(x)},0\bigr)-
   \bigl(\overline{\bar T_s(x)},0\bigr)=
   \epsilon_2\bigl(\bar T_s(x)\bigr).
\end{split}
\]
(Note that $T_s=L_s+2R_s$, so $\bar T_s=-\bigl(R_s+2L_s\bigr)$.)

From Lemma \ref{le:lrt} we know that $\bigl(-R_s,L_s+R_s,-L_s\bigr)$
and $\bigl(-L_s,-R_s,L_s+R_s\bigr)$ belong to $\lrt(A,\bar\ )$, and
so does their difference $\bigl(L_s-R_s,T_s,\bar T_s\bigr)$. Also,
since any element $(d_0,d_1,d_2)$ in $\lrt(A,\bar\ )$ is determined
by the pair $(d_1,d_2)$, it follows that
\[
\psi(T_s)=\bigl(L_s-R_s,T_s,\bar T_s\bigr).
\]

Moreover,
\[
\begin{split}
[\alpha(0,s)+\alpha^{-1}(0,s)\tilde\ ,\epsilon_1(x)]
 &=[\alpha(0,s)+\alpha^{-1}(0,s)\tilde\ ,
   (x,0)+\alpha^{-1}(x,0)\tilde\ ]\\
 &=[(0,s),(x,0)\tilde\ ]-\alpha^{-1}[(x,0),(0,s)\tilde\ ]\\
 &=(sx,0)+\alpha^{-1}(sx,0)\tilde\ =\epsilon_1(sx),\\[2pt]
[\alpha(0,s)+\alpha^{-1}(0,s)\tilde\ ,\epsilon_2(x)]
 &=[\alpha(0,s)+\alpha^{-1}(0,s)\tilde\ ,
    \alpha(\bar x,0)-(\bar x,0)\tilde\ ]\\
 &=-\alpha[(0,s),(\bar x,0)\tilde\ ]-
   [(\bar x,0),(0,s)\tilde\ ]\\
 &=-\alpha(s\bar x,0)+(s\bar x,0)\tilde\ =\epsilon_2(xs),
\end{split}
\]
whence we conclude that
\[
\psi\bigl(\alpha(0,s)+\alpha^{-1}(0,s)\tilde\ \bigr)=
 -\bigl(-(L_s+R_s),L_s,R_s\bigr).
\]
Therefore, for any $s_0,s_1,s_2\in S$ with $s_0+s_1+s_2=0$,
\begin{multline*}
\bigl(L_{s_1}-R_{s_2},L_{s_2}-R_{s_0},L_{s_0}-R_{s_1}\bigr)\\
 =\frac{1}{2}\Bigl(\psi\bigl(T_{s_1+s_2}\bigr)
   -\psi\bigl(\alpha(0,s_1-s_2)+\alpha^{-1}(0,s_1-s_2)\tilde\
   \bigr)\Bigr),
\end{multline*}
so that
\[
\psi\Bigl(T_S\oplus \{\alpha(0,s)+\alpha^{-1}(0,s)\tilde\ :s\in
S\}\Bigr)=\calT_S,
\]
and hence $\psi$ gives an isomorphism onto
$\frl(A,-\frd)=\frd^{<3>}\oplus\calT_S$.
\end{proof}

\medskip

Let us recall now Allison and Faulkner's construction of the Lie
algebra $\calK(A,\bar\ ,\gamma,\frv)$, for any structurable algebra
$(A,\bar\ )$, $\gamma=(\gamma_1,\gamma_2,\gamma_3)\in F^3$ and
subalgebra $\frv$ of $\lrt(A,\bar\ )$ containing $\inlrt(A,\bar\ )$.
As a vector space, $\calK(A,\bar\ ,\gamma,\frv)=\frv\oplus
A[12]\oplus A[23]\oplus A[31]$, where $A[ij]$ is a copy of $A$ with
$a[ij]=-\gamma_i\gamma_j^{-1}\bar a[ji]$ for any $1\leq i\ne j\leq
3$, and the multiplication is obtained by extending the bracket in
$\frv$ by setting for any $a,b\in A$
\[
\begin{split}
\bigl[a[ij],b[jk]\bigr]&=ab[ik],\\
\bigl[T,a[ij]\bigr]&=T_k(a)[ij],\ \text{for $T\in\frv$}\\
\bigl[a[ij],b[ij]\bigr]&=\gamma_i\gamma_j^{-1}T,
\end{split}
\]
where $(i,j,k)$ is a cyclic permutation of $(1,2,3)$ and the $T$ in
the last row is $T=(T_1,T_2,T_3)$ with
\[
\begin{split}
T_i&=L_{\bar b}L_a-L_{\bar a}L_b,\\
T_j&=R_{\bar b}R_a-R_{\bar a}R_b,\\
T_k&=R_{(\bar ab-\bar ba)}+L_bL_{\bar a}-L_aL_{\bar b}.
\end{split}
\]
(Compare with Example \ref{ex:allison}.)

\smallskip

\begin{proposition}\label{pr:Psi}
Let $(A,-)$ be a structurable algebra and let $\frd$ be a subalgebra
of $\der(A,\bar\ )$ containing the inner derivations. Then
$\calK(A,\bar\ ,\frd)$ is isomorphic to $\calK\bigl(A,\bar\
,\gamma,\frl(A,\bar\ ,\frd)\bigr)$, with $\gamma=(1,-1,2\alpha)$.
\end{proposition}
\begin{proof}
Consider the following elements in $\calK(A,\bar\ ,\frd)$:
\[
\tilde\epsilon_0(x)=\epsilon_0(x),\quad
\tilde\epsilon_1(x)=\frac{1}{2}\epsilon_1(x),\quad
\tilde\epsilon_2(x)=-\epsilon_2(x),
\]
for any $x\in A$. Then \eqref{eq:epsilones} becomes
\begin{align*}
[\tilde\epsilon_0(x),\tilde\epsilon_1(y)]&
      =-\frac{1}{2\alpha}\tilde\epsilon_2(\overline{xy}),\\
[\tilde\epsilon_1(x),\tilde\epsilon_2(y)]&=\tilde\epsilon_0(\overline{xy}),\\
[\tilde\epsilon_2(x),\tilde\epsilon_0(y)]&
    =2\alpha\tilde\epsilon_1(\overline{xy}).
\end{align*}
Note that for $\gamma=(1,-1,2\alpha)$,
$-\frac{1}{2\alpha}=-\gamma_1\gamma_3^{-1}$,
$1=-\gamma_2\gamma_1^{-1}$ and $2\alpha=-\gamma_3\gamma_2^{-1}$.
Now, the isomorphism $\psi:\calK_{(\bar 0,\bar 0)}\rightarrow
\frl(A,\bar\ ,\frd)$, can be extended to an isomorphism
$\Psi:\calK(A,\bar\ ,\frd)\rightarrow \calK\bigl(A,\bar\
,\gamma,\frl(A,\bar\ ,\frd)\bigr)$ by means of
\[
\Psi\bigl(\tilde\epsilon_0(x)\bigr)=x[12],\quad
\Psi\bigl(\tilde\epsilon_1(x)\bigr)=x[23],\quad
\Psi\bigl(\tilde\epsilon_2(x)\bigr)=x[31],
\]
while, for any $p\in \calK_{(\bar 0,\bar 0)}$,
$\Psi(p)=(d_1,d_2,d_0)$ if $\psi(p)=(d_0,d_1,d_2)$.  The only
difficulty in proving that $\Psi$ is an isomorphism lies in proving
that $\Psi\bigl[\tilde\epsilon_i(x),\tilde\epsilon_i(y)\bigr]
=\bigl[\Psi\bigl(\tilde\epsilon_i(x)\bigr),
\Psi\bigl(\tilde\epsilon_i(y)\bigr)\bigr]$ for any $i\in\bZ_3$, and
$x,y\in A$. But the action of both sides on
$\Psi\bigl(\tilde\epsilon_{i+1}(A)\oplus\tilde\epsilon_{i+2}(A)\bigr)$
coincide, while any element in $\lrt(A,\bar\ )$ is determined by its
action on two of the direct summands in $A[12]\oplus A[23]\oplus
A[31]=\Psi\bigl(\tilde\epsilon_0(A)\oplus\tilde\epsilon_1(A)\oplus
\tilde\epsilon_2(A)\bigr)$.
\end{proof}

\medskip

\begin{corollary}
Let $(A,-)$ be a structurable algebra over a field $F$ satisfying
that $-1\in F^2$, and let $\frd$ be a subalgebra of $\der(A,\bar\ )$
containing the inner derivations. Then, $\Aut\bigl(\calK(A,\bar\
,\frd)\bigr)$ contains a subgroup isomorphic to the symmetric group
of degree $4$.
\end{corollary}
\begin{proof}
Take $\alpha=2$ above and define
\[
\iota_0(x)=\sqrt{-1}\epsilon_0(x),\quad
\iota_1(x)=\sqrt{-1}\epsilon_1(x),\quad
\iota_2(x)=\frac{1}{2}\epsilon_2(x),
\]
for any $x\in A$. Then \eqref{eq:epsilones} becomes
\[
[\iota_i(x),\iota_{i+1}(y)]=\iota_{i+2}(\overline{xy})
\]
for any $x,y\in A$ and $i\in\bZ_3$. With the arguments of the last
proof, it is readily seen that $\calK(A,\bar\ ,\frd)$ is isomorphic
to $\calK\bigl(A,\bar\ ,\gamma,\frl(A,\bar\ ,\frd)\bigr)$, where
$\gamma=(-1,-1,-1)$. Now, besides the automorphisms $\tau_1,\tau_2$
used to obtain the grading over $\bZ_2\times\bZ_2$, there appears
the order $3$ automorphism $\varphi$ such that
\[
\left\{\begin{aligned}
&\varphi\bigl(\iota_i(x)\bigr)=\iota_{i+1}(x),\\
&\varphi\bigl(\psi^{-1}(d_0,d_1,d_2)\bigr)=\psi^{-1}(d_2,d_0,d_1)
\end{aligned}\right.
\]
and the order two automorphism $\tau$ such that
\[
\left\{\begin{aligned}
 &\tau\bigl(\iota_0(x)\bigr)=-\iota_0(\bar x),\\
 &\tau\bigl(\iota_1(x)\bigr)=-\iota_2(\bar x),\\
 &\tau\bigl(\iota_2(x)\bigr)=-\iota_1(\bar x),\\
 &\tau\bigl(\psi^{-1}(d_0,d_1,d_2)\bigr)=\psi^{-1}(\bar d_0,\bar
 d_2,\bar d_1),
\end{aligned}\right.
\]
with $x\in A$, $i\in \bZ_3$ and $(d_0,d_1,d_2)\in \frl(A,\bar\
,\frd)$. These automorphisms generate a subgroup of the automorphism
group isomorphic to $S_4$.
\end{proof}

\medskip

It must be remarked that Allison proved in \cite[Theorem 2.2 and
Section 4]{All91} that if $-\gamma_1\gamma_2^{-1}\in (F^\times)^2$,
then $\calK\bigl(A,\bar\ ,\gamma,\frl(A,\bar\ ,\frd)\bigr)$ is
isomorphic to $\calK(A,\bar\ ,\frd)$. In particular, $\calK(A,\bar\
,\frd)$ is isomorphic to $\calK\bigl(A,\bar\ ,(1,-1,1),\frl(A,\bar\
,\frd)\bigr)$.


\providecommand{\bysame}{\leavevmode\hbox
to3em{\hrulefill}\thinspace}
\providecommand{\MR}{\relax\ifhmode\unskip\space\fi MR }
\providecommand{\MRhref}[2]{%
  \href{http://www.ams.org/mathscinet-getitem?mr=#1}{#2}
} \providecommand{\href}[2]{#2}

\end{document}